\documentclass[10pt]{article}

\usepackage{graphicx,hyperref,slashbox}
\usepackage{amsmath}
\usepackage[version=4]{mhchem}
\usepackage{siunitx}
\usepackage{longtable,tabularx}
\usepackage{pstricks}
\usepackage{amsfonts}
\usepackage{latexsym}
\usepackage{enumerate}
\usepackage{epsfig}
\usepackage{subfig}
\usepackage{framed,fancybox}
\usepackage{bbm}
\usepackage{multirow}
\usepackage[T1]{fontenc}
\usepackage{hhline}

\usepackage{threeparttable}
\usepackage{rotating}
\usepackage{amssymb,amscd}
\makeatletter
\def\myVCENTER#1{\vcenter{\hbox{$\m@th#1$}}}
\makeatother

\long\def\symbolfootnote[#1]#2{\begingroup\def\thefootnote{\fnsymbol{footnote}}\footnote[#1]{#2}\endgroup}

\definecolor{shadecolor}{gray}{0.99}
{\endMakeFramed}

\newenvironment{shadedframe}{%
 \MakeFramed {\FrameRestore}}
{\endMakeFramed}

\definecolor{shadecolor}{gray}{0.99}
\long\def\symbolfootnote[#1]#2{\begingroup\def\thefootnote{\fnsymbol{footnote}}\footnote[#1]{#2}\endgroup}
\def\qed{\hfill{$\vcenter{\hrule height1pt \hbox{\vrule width1pt height5pt
    \kern5pt \vrule width1pt} \hrule height1pt}$} \medskip}
\newcommand{\m}[1]{{\bf{#1}}}
\newcommand{\g}[1]{\boldsymbol #1}
\newcommand{\bb}[1]{\mathbb #1}

\newcommand{\C}[1]{{\cal {#1}}}

\title{\bf A Warm Start Method for Solving \\ Chance Constrained Optimal Control Problems}

\author{Rachel E.~Keil\thanks{Ph.D.~Candidate, Department of Mechanical and Aerospace Engineering, University of Florida, Gainesville, FL 32611-6250. E-mail: rekeil@ufl.edu} \\ Mrinal Kumar \thanks{Associate Professor, Department of Mechanical and Aerospace Engineering. The Ohio State University, Columbus, OH 43210. AIAA Senior Member. E-mail: kumar.672@osu.edu} \\ Anil V.~Rao\thanks{Professor, Erich Farber Faculty Fellow and University Term Professor, Department of Mechanical and Aerospace Engineering, University of Florida, Gainesville, FL 32611-6250.  E-mail: anilvrao@ufl.edu.  Corresponding Author.}}
\date{}
\begin{document}
\maketitle
\thispagestyle{empty}

\begin{abstract}
A warm start method is developed for efficiently solving complex chance constrained optimal control problems.  The warm start method addresses the computational challenges of solving chance constrained optimal control problems using biased kernel density estimators and Legendre-Gauss-Radau collocation with an $hp$ adaptive mesh refinement method.  To address the computational challenges, the warm start method improves both the starting point for the chance constrained optimal control problem, as well as the efficiency of cycling through mesh refinement iterations.  The improvement is accomplished by tuning a parameter of the kernel density estimator, as well as implementing a kernel switch as part of the solution process.  Additionally, the number of samples for the biased kernel density estimator is set to incrementally increase through a series of mesh refinement iterations.  Thus, the warm start method is a combination of tuning a parameter, a kernel switch, and an incremental increase in sample size.  This warm start method is successfully applied to solve two challenging chance constrained optimal control problems in a computationally efficient manner using biased kernel density estimators and Legendre-Gauss-Radau collocation.    
\end{abstract}

\newpage

\renewcommand{\baselinestretch}{2}
\normalsize\normalfont 

\section*{Nomenclature}

{\renewcommand\arraystretch{1.0}
\noindent\begin{longtable*}{@{}l @{\quad=\quad} l@{}}
$(\cdot)^*$ & symbol for optimal solution \\
$t$ & time \\
$\tau $ & time interval transformation \\
$J$ & cost function \\
$\C{M}$ & Mayer cost \\
$\C{L}$ & Lagrangian \\
$\m{y}$ & state \\
$\m{Y}$ & state approximation\\
$\m{u}$ & control \\
$\m{U}$ & control approximation\\
$\m{c}$ & path constraint \\
$\m{b}$ & boundary constraint \\
$\m{g}$ & chance constraint \\
$K$ & number of mesh intervals \\
$N_k$ & number of LGR collocation points in mesh interval $k$ \\
$\ell_i$ & $i^{th}$ Lagrange polynomial \\
$t$ & time \\
$\m{D}$ & Differentiation matrix \\ 
$P_N$ & $N^{th}$--degree Legendre Polynomial \\
$\g{\lambda}$ & costate \\
$\g{\Lambda}$ & defect constraint Lagrange multiplier \\
$P$ & probability  \\
$\Omega$ & sample space for a random event \\
$\bb{E}$ & expectation (mean) function \\
$1_{(\cdot)}$ & indicator function \\
$\boldsymbol{\xi}$ & random vector \\
$\psi$ & random variable \\
$ \mathbf{q} $ & event boundary \\
$\epsilon_{()}$ & risk violation parameter \\
$k_{()}$ & kernel \\
$K_{()}$ & integrated kernel function \\
$h$ & bandwidth \\
$\phi$ & user defined parameter \\
$B$ & bias of kernel \\
$T$ & execution time for $\mathbb{GPOPS-II}$ \\
$C$ & number of convergent runs \\
$H$ & number of runs converging to higher cost solution \\

\end{longtable*}}

\section{Introduction}
Optimal control problems arise frequently in a variety of engineering and non-engineering disciplines.  The goal of an optimal control problem is to determine the state and control of a controlled dynamical system that optimize a given performance index while being subject to path constraints, event constraints, and boundary conditions~\cite{Betts3}.  Optimal control problems are either deterministic or stochastic.  A deterministic optimal control problem is an optimal control problem that contains no uncertainty, while a stochastic optimal control problem is one that contains uncertainty.  Forms of uncertainty include measurement error, process noise, model error, and uncertainty in the constraints.  Examples where the constraints contain uncertainty include fuzzy boundary keep out zone path constraints~\cite{Keil1}, variable control limitation path constraints~\cite{Kumar1}, and event constraints with variations in the state~\cite{Kumar1,Caillau}.  Constraints with uncertainty are often modeled as chance constraints, and optimal control problems subject to chance constraints are called chance constrained optimal control problems (CCOCPs).  

Due to the probabilistic form of the chance constraints, most CCOCPs must be solved numerically.  Numerical methods for solving optimal control problems have been, however, developed primarily for solving deterministic optimal control problems.  As a result, methods for transforming the CCOCP to a deterministic optimal control problem have been developed.  Many of these methods focus on transforming the chance constraints to deterministic constraints in a manner that retains the key stochastic properties of the original chance constraint.  Such methods include the methods of Refs.~\cite{blackmore10,Blackmore1,ono10,okamoto19,hokayem13} that are applicable to linear chance constraints.  The methods of Refs.~\cite{Pinter,muhlpfordt18,Nemirovski} are applicable to chance constraints when certain information about the chance constraint is available.  When information about the chance constraints is not available, the methods of Refs.~\cite{Kumar1,Caillau,Pagnoncelli1,ono15,Calafiore1,Calafiore2,Campi1,Chai,Ahmed,Calfa,Keil2} are applicable.  
  
Recently, Ref.~\cite{Keil2} has developed a new method for transforming chance constraints to deterministic constraints using biased kernel density estimators (KDEs).  An advantage of the method developed in Ref.~\cite{Keil2} is that the deterministic constraint is not overly conservative relative to the chance constraint and does not violate the boundary of the chance constraint.  In addition, the method developed in Ref.~\cite{Keil2} has a key feature that it is formulated using an adaptive Gaussian quadrature orthogonal collocation method \cite{Benson2,Rao8,Garg1,Garg2,Patterson2015} known as LGR collocation.  By combining biased KDEs with LGR collocation as performed in Ref.~\cite{Keil2}, it is possible to take advantage of several properties of Gaussian quadrature collocation.  First, using Gaussian quadrature collocation, the constraints are evaluated independently at each collocation point.  In addition, Gaussian quadrature collocation provides high-accuracy solutions along with exponential convergence for smooth optimal control problems \cite{Garg1,Garg2,Patterson2015,HagerHouRao15a,HagerHouRao16a,HagerLiuMohapatraWangRao19,DuChenHager2019}.

While Ref.~\cite{Keil2} provides a method for transforming CCOCPS to deterministic optimal control problems using biased KDEs, it does not provide a computationally reliable and efficient method for solving the resulting optimization problem.  In fact, it is shown in this paper that, using a naive approach, solving the nonlinear programming problem (NLP) that arises from the approach of Ref.~\cite{Keil2} produces different results depending upon the manner in which the problem is initialized.  Moreover, even when a solution can be obtained, it is shown that unnecessary computational effort is required.  As a result, it is important in practical applications to to develop a computational approach that can be used with the method of Ref.~\cite{Keil2} that simultaneously leads to a tractable optimization problem and enables the optimization problem to be solved efficiently.  Situations that could benefit from the approach developed in this paper include rapid (that is, on short notice) or real time trajectory optimization.  

This paper describes a new computational approach for solving CCOCPs using biased KDEs together with LGR collocation.  The approach developed in this paper is called a {\em warm start method} because it improves the initial guess provided to the NLP solver for the transcribed CCOCP.  Moreover, because the transcribed CCOCP is solved using collocation together with mesh refinement, this warm start approach is used to generate an initial guess for the NLP solver on each mesh.  The key benefit of this approach is that it improves the reliability and efficiency of solving CCOCPs using biased KDEs together with LGR collocation.  

The warm start method developed in this paper has three major components.  The first component tunes a parameter of the biased KDE in order to improve the starting point for the NLP.  The second component is a kernel switching procedure that allows changing kernels which ensures that the starting kernel leads to a tractable optimization problem, while maintaining the ability to obtain results for other kernels.  The third component is a procedure that incrementally increases the number of samples required for use with the biased KDEs.

The key contribution of this research is to develop a novel method to reliably and efficiently solve the NLP that arises from transforming a CCOCP using biased KDEs together with LGR collocation.  The goal is to provide researchers with an approach that is tractable for solving increasingly complex CCOCPs.  Using a systematic formulation for the method to solve two example CCOCPs, significant improvements are found using the approach developed in this paper.  

This paper is organized as follows.  Section~\ref{sect:review} provides a brief review of biased KDEs and LGR collocation.  In Section~\ref{sect:Jorris2Dnaive}, a complex CCOCP is solved using biased KDEs and LGR collocation without a warm start.  In Section~\ref{sect:tech}, the warm start method is developed.  Section~\ref{sect:discusstech} provides the results of solving the complex CCOCP with the warm start method.   In Section~\ref{sect:examples}, a second complex CCOCP is solved using the warm start method.  Sections~\ref{sect:conclude} and~\ref{sect:discussion} provide a discussion and some conclusions, respectively.  

\section{Chance Constrained Optimal Control\label{sect:review}}

In this section, biased kernel density estimators (KDEs) are combined with Legendre-Gauss-Radau (LGR) collocation to transform a chance constrained optimal control problem (CCOCP) to a nonlinear programming problem (NLP).   First, Section~\ref{sect:discCCOCP} describes a general continuous CCOCP.  Section~\ref{sect:LGRcolloc} then describes LGR collocation.  Finally, in Section~\ref{sect:biasKDE}, biased KDEs are applied to transform the chance constraints to deterministic constraints, where the deterministic constraints retain the main probability properties of the chance constraint.

\subsection{General Chance Constrained Optimal Control Problem}\label{sect:discCCOCP}
Consider the following general continuous CCOCP.  Determine the state
$\m{y}(\tau)\in\bb{R}^{n_y}$ and the control $\m{u}(\tau)\in\bb{R}^{n_u}$ on the domain $\tau \in [-1, +1]$,  the initial time, $t_0$, and the terminal time $t_f$ that minimize the cost functional

\begin{subequations}
\begin{equation}\label{bolza-cost-s}
  \C{J} =\C{M}(\m{y}(-1),t_0,\m{y}(+1),t_f)  + \frac{t_f-t_0}{2}\int_{-1}^{+1} \C{L}(\m{y}(\tau),\m{u}(\tau), t(\tau, t_0, t_f))\, d\tau,
\end{equation}
subject to the dynamic constraints
\begin{equation}\label{bolza-dyn-s}
  \frac{d\m{y}}{d\tau} -
\frac{t_f-t_0}{2}\m{a}(\m{y}(\tau),\m{u}(\tau), t(\tau, t_0, t_f) )=\m{0}, 
\end{equation}
the inequality path constraints
\begin{equation}\label{bolza-path-s}
\m{c}_{\min} \leq \m{c}(\m{y}(\tau),\m{u}(\tau), t(\tau, t_0, t_f) )\leq \m{c}_{\max},
\end{equation}
the boundary conditions
\begin{equation}\label{bolza-bc-s}
  \m{b}_{\min} \leq \m{b}(\m{y}(-1),t_0,\m{y}(+1),t_f) \leq \m{b}_{\max},
\end{equation}
and the chance constraint
\begin{equation}\label{bolza-pathcc-s}
P( \m{F} (\m{y}(\tau),\m{u}(\tau),t(\tau, t_0, t_f);\g{\xi}) \geq \m{q}) \geq 1- \epsilon.
\end{equation}
\end{subequations}
The random vector $\g{\xi}$ is supported on set $\Omega \subseteq \bb{R}^d$.  The function $ \m{F}(\cdot) \geq \m{q}$ is an event in the probability space $P(\cdot)$ where $\m{q} \subset \bb{R}^{n_g}$ is the boundary of the event and $\epsilon$ is the risk violation parameter.  It is noted that the path, event and dynamic constraints can all be in the form of a chance constraint.  Additionally, it is noted that the time interval $\tau\in[-1,+1]$ can be transformed
to the time interval $t\in[t_0,t_f]$ via the affine transformation
\begin{equation}\label{tau-to-t}
  t \equiv t(\tau,t_0,t_f) = \frac{t_f-t_0}{2}\tau + \frac{t_f+t_0}{2}.  
\end{equation}

The continuous CCOCP of Eqs.~\eqref{bolza-cost-s}-\eqref{bolza-pathcc-s} must be transformed to a form that is solvable using numerical methods.  For application with numerical methods, the CCOCP is discretized on the domain $\tau\in[-1,+1]$ which is partitioned into a {\em mesh} consisting of $K$  {\em mesh intervals} $\C{S}_k=[T_{k-1},T_k],\; k=1,\ldots,K$, where $-1 = T_0 < T_1 < \ldots < T_K = +1$.  The mesh intervals have the property that $\displaystyle \cup_{k=1}^{K} \C{S}_k=[-1,+1]$.  Let $\m{y}^{(k)}(\tau)$ and $\m{u}^{(k)}(\tau)$ be the state and control in $\C{S}_k$.  Using the transformation given in Eq.~\eqref{tau-to-t}, the CCOCP of Eqs.~\eqref{bolza-cost-s}-\eqref{bolza-pathcc-s} can then be rewritten as follows.  Minimize the cost functional

\begin{subequations}
\begin{equation}\label{bolza-cost-segmented}
    \C{J} = \C{M}(\m{y}^{(1)}(-1),t_0,\m{y}^{(K)}(+1),t_f)  + \frac{t_f-t_0}{2}\sum_{k=1}^K \int_{T_{k-1}}^{T_k} \C{L}(\m{y}^{(k)}(\tau),\m{u}^{(k)}(\tau),t)\, d\tau,
\end{equation}
subject to the dynamic constraints
\begin{equation}\label{bolza-dyn-segmented}
  \displaystyle\frac{d\m{y}^{(k)}(\tau)}{d\tau} - \frac{t_f-t_0}{2}\m{a}(\m{y}^{(k)}(\tau),\m{u}^{(k)}(\tau), t)=\m{0}, \quad  (k=1,\ldots,K),
\end{equation}
the path constraints
\begin{equation}\label{bolza-path-segmented}
  \m{c}_{\min} \leq \m{c}(\m{y}^{(k)}(\tau),\m{u}^{(k)}(\tau), t) \leq \m{c}_{\max},\quad  (k=1,\ldots,K),
\end{equation}
the boundary conditions
\begin{equation}\label{bolza-bc-segmented}
\m{b}_{\min} \leq \m{b}(\m{y}^{(1)}(-1),t_0,\m{y}^{(K)}(+1),t_f) \leq  \m{b}_{\max}, 
\end{equation}
and the chance constraint
\begin{equation}\label{bolza-pathcc-segmented}
P(\m{F} (\m{y}^{(k)}(\tau),\m{u}^{(k)}(\tau),t;\g{\xi}) \geq \m{q}) \geq 1- \epsilon.
\end{equation}
Because the state must be continuous, the following condition 
\begin{equation}\label{eq:contin}
  \m{y}^{(k)}(T_k)=\m{y}^{(k+1)}(T_k),\;(k=1,\ldots,K-1).
\end{equation}
\end{subequations}

\subsection{Legendre-Gauss-Radau Collocation}\label{sect:LGRcolloc}
The form of discretization that will be applied to the CCOCP in Section~\ref{sect:discCCOCP} is collocation at
LGR points~\cite{Garg1,Garg2,Patterson2015}.  In the LGR 
collocation method, the state of the continuous CCOCP is approximated in $\C{S}_k,\;k\in[1,\ldots,K]$, as 
\begin{equation}\label{state-approximation-LGR}
\begin{split}
\m{y}^{(k)}(\tau)  \approx \m{Y}^{(k)}(\tau) & = \sum_{j=1}^{N_k+1}
\m{Y}_{j}^{(k)} \ell_{j}^{(k)}(\tau),\\ \ell_{j}^{(k)}(\tau) & = \prod_{\stackrel{l=1}{l\neq j}}^{N_k+1}\frac{\tau-\tau_{l}^{(k)}}{\tau_{j}^{(k)}-\tau_{l}^{(k)}}, 
\end{split}
\end{equation}  
where $\tau\in[-1,+1]$, $\ell_{j}^{(k)}(\tau),$ $j=1,\ldots,N_k+1$, is a
basis of Lagrange polynomials,
$\left(\tau_1^{(k)},\ldots,\tau_{N_k}^{(k)}\right)$ are the 
LGR collocation points in $\C{S}_k =$ $[T_{k-1},T_k)$, and 
$\tau_{N_k+1}^{(k)}=T_k$ is a noncollocated point.  Differentiating
$\m{Y}^{(k)}(\tau)$ in Eq.~(\ref{state-approximation-LGR}) with
respect to $\tau$ gives
\begin{equation}\label{diff-state-approximation-LGR}
  \frac{d\m{Y}^{(k)}(\tau)}{d\tau} = \sum_{j=1}^{N_k+1}\m{Y}_{j}^{(k)}\frac{d\ell_j^{(k)}(\tau)}{d\tau}.
\end{equation}
Defining $t_i^{(k)}=t(\tau_i^{(k)},t_0,t_f)$ using
Eq.~\eqref{tau-to-t}, the dynamics are then approximated at the $N_k$
LGR points in mesh interval $k\in[1,\ldots,K]$ as
\begin{equation}\label{collocation-LGR}
    \sum_{j=1}^{N_k+1}D_{ij}^{(k)} \m{Y}_j^{(k)} - \frac{t_f-t_0}{2}\m{a}(\m{Y}_i^{(k)},\m{U}_i^{(k)},t_i^{(k)})=\m{0}, \  (i=1,\ldots,N_k),
\end{equation}
where $D_{ij}^{(k)} = d\ell_j^{(k)}(\tau_i^{(k)})/d\tau,\;(i=1,\ldots,N_k),\;(j=1,\ldots,N_k+1)$ are the elements of the $N_k\times (N_k+1)$ {\em Legendre-Gauss-Radau differentiation matrix}~\cite{Garg1} in mesh interval
$\C{S}_k,\;k\in[1,\ldots,K]$.  The LGR discretization then leads to
the following resulting form of the discretized CCOCP.  Minimize 
\begin{equation}\label{cost-LGR}
    \C{J}  \approx \C{M}(\m{Y}_{1}^{(1)},t_0,\m{Y}_{N_K+1}^{(K)},t_f) +   \sum_{k=1}^{K} \sum_{j=1}^{N_k}  \frac{t_f-t_0}{2} 
  w_{j}^{(k)} \C{L}(\m{Y}_{j}^{(k)},\m{U}_{j}^{(k)},t_j^{(k)}),
\end{equation}
subject to the collocation constraints of Eq.~\eqref{collocation-LGR}
and the constraints
\begin{gather}\label{eq:differential-collocation-conditions-LGR}
  \m{c}_{\min} \leq \m{c}(\m{Y}_{i}^{(k)},\m{U}_{i}^{(k)},t_i^{(k)}) \leq \m{c}_{\max},\; (i=1,\ldots,N_k),\\
 \m{b}_{\min} \leq \m{b}(\m{Y}_{1}^{(1)},t_0,\m{Y}_{N_K+1}^{(K)},t_f)  \leq \m{b}_{\max},  \\
 P( \m{F} (\m{Y}_{i}^{(k)},\m{U}_{i}^{(k)},t_i^{(k)};\g{\xi}) \geq \m{q}) \geq 1- \epsilon ,\; (i=1,\ldots,N_k) , \label{cc-constraint}\\
\m{Y}_{N_k+1}^{(k)} = \m{Y}_1^{(k+1)} , \quad (k=1,\ldots,K-1),  \label{continuity-constraint}
\end{gather}
where $N = \sum_{k=1}^{K} N_k$ is the total number of LGR points and 
Eq.~\eqref{continuity-constraint} is the continuity condition on
the state that is enforced at the interior mesh points
$(T_1,\ldots,T_{K-1})$ by treating $\m{Y}_{N_k+1}^{(k)}$ and
$\m{Y}_1^{(k+1)}$ as the same variable.  In order for Eqs.~\eqref{collocation-LGR}-~\eqref{continuity-constraint} to be an NLP, the chance constraint is transformed to a deterministic constraint in the next section.

\subsection{Biased Kernel Density Estimators}\label{sect:biasKDE}
In this section, the chance constraint of Eq.~\eqref{cc-constraint} is transformed to a deterministic constraint using biased KDEs~\cite{Keil2}.  In order to transform the chance constraint, first, it is noted that the function $\m{F}(\cdot)$ of Eq.~\eqref{cc-constraint} is itself a random vector $\g{\psi}$ whose associated probabilistic properties are unknown.  Consequently, the constraint of Eq.~\eqref{cc-constraint} can be redefined as: 
\begin{equation}\label{eq:CCexample}
P(\g{\psi} \geq \m{q}) \geq 1-\epsilon.
\end{equation}
Because $\g{\psi}$ is a random vector, Eq.~\eqref{eq:CCexample} is a joint chance constraint.  Using Boole's inequality together with the approach of Refs.~\cite{Nemirovski,Blackmore1,Kumar1}, the chance constraint given in Eq.~\eqref{eq:CCexample} can be redefined in terms of the following set of conservative constraints (see Refs.~\cite{Blackmore1} and~\cite{Kumar1} for the proof):
\begin{equation}\label{eq:CC_scalar}
\begin{array}{c}
P( \psi \geq q_m  ) \geq  1-\epsilon_m, \\
\sum\limits_{m =1}^{n_g} \epsilon_m \leq \epsilon,
\end{array}
\end{equation} 
or, equivalently
\begin{equation}\label{eq:CC_scalar_comp}
\begin{array}{c}
P( \psi < q_m  ) \leq  \epsilon_m, \\
\sum\limits_{m =1}^{n_g} \epsilon_m \leq \epsilon,
\end{array}
\end{equation}
where $m\in [1,\dots,n_g]$ is the index corresponding to the $m$th component of the event and $\psi$ is the $m$th variable of random vector $\g{\psi}$.

As a result of the chance constraint now being in the scalar form of Eq.~\eqref{eq:CC_scalar}, this chance constraint can be transformed to a deterministic constraint using biased KDEs.  First, for this transformation, the kernel $k(\cdot) $ of the KDE is integrated to obtain the following integrated kernel function $K(\cdot)$:
\begin{equation}\label{eq:kernCDF}
\begin{array}{c}
K(\eta_j) = \int_{-\infty}^{\eta_j}k(v_j) d v_j, \\
\eta_j = \frac{q_m-\psi_j}{h},
\end{array}
\end{equation}
where $h$ is the bandwidth and $j = 1,\dots,N$ is the number of samples of the random variable $\psi$.  It is noted that the samples of $\psi$ are obtained by sampling $\g{\xi}$.  Next, the following relation between the chance constraint of Eq.~\eqref{eq:CC_scalar} and the integrated kernel function $K(\cdot)$ is defined:
\begin{equation}\label{eq:KDEbias}
\frac{1}{N} \sum_{j = 1}^N K_B \left( \eta_j \right) \leq P(\psi < q_m) \leq \epsilon_m,
\end{equation}
where the subscript $B$ on the integrated kernel function indicates that the kernel function has been biased by amount $B(h)$ and the left-hand side of Eq.~\eqref{eq:KDEbias} is the biased KDE.
The bias $B(h)$ is chosen such that the biased integrated kernel function $K_B(\cdot)$ satisfies the following inequality:
\begin{align}\label{eq:geninequality}
1_{[0,+\infty)} \left( \nu \right) \leq K_B(\nu), \ \forall \nu,\\
\nu = \frac{\psi-q_m}{h}, \ \textrm{with}~h >0.
\end{align} 
The relation of Eq.~\eqref{eq:geninequality} is the first requirement from Ref.~\cite{Keil2} for the relation between the biased KDE and the chance constraint from Eq.~\eqref{eq:KDEbias} to hold.  The second requirement is that the number of MCMC samples $N$ reaches a value $N_c$ that is sufficiently large~\cite{Kumar1,Keil2,Kumar3,Kumar2} in order to accurately approximate the characteristics of the distribution of the random vector $\boldsymbol{\xi}$ \cite{MCMCMethods}.  If these samples are available, the following expression  
\begin{equation}\label{eq:inequaltosat}
\lim_{N \to N_c} \frac{1}{N} \sum\limits_{j = 1}^N K_B (\nu_j) \ = \bb{E} [ K_B (\nu)] \geq 1-\epsilon_m, 
\end{equation}
converges to the expectation $\bb{E} [K_B(\nu)]$, where the expectation $\bb{E} [ K_B (\nu)]$ exists for $h>0$ for the nonempty compact set
\begin{equation}\label{eq:setsforproof}
1-\epsilon_m \leq \bb{E} [ K_B (\nu)].
\end{equation}
If the first and second requirements from, respectively, Eq.~\eqref{eq:geninequality} and Eq.~\eqref{eq:inequaltosat} are satisfied, the chance constraint of Eq.~\eqref{cc-constraint} can be transformed to the following set of deterministic constraints:
\begin{equation}\label{eq:KDEfinalform}
\begin{array}{c}
\frac{1}{N} \sum\limits_{j = 1}^N K_B \left( \eta_j \right) \leq \epsilon_m, \\
\sum\limits_{m =1}^{n_g} \epsilon_m \leq \epsilon.
\end{array}
\end{equation}
By replacing the chance constraint from Eq.~\eqref{cc-constraint} by the deterministic constraints from Eq.~\eqref{eq:KDEfinalform}, the system of equations from Eqs.~\eqref{collocation-LGR}--\eqref{continuity-constraint} is now a NLP that can be solved available software such as SNOPT~\cite{Gill1,Gill2}, IPOPT~\cite{Biegler2}, and KNITRO~\cite{Byrd1}. 

\section{Motivation for Warm Start Method\label{sect:Jorris2Dnaive}}

This section provides motivation for the warm start method developed in Section~\ref{sect:tech}.  This motivation is furnished via a complex example CCOCP that is solved using biased KDEs and LGR collocation without a warm start method.  In Section~\ref{sect:Jorris2D}, the example is presented.  Section~\ref{sect:guesseschoic} provides the initialization.  Next, Section~\ref{sect:setup} describes the setup for the optimal control software $\mathbb{GPOPS-II}$.  Finally, Section~\ref{sect:Jorris2Ddiscuss} provides the results of solving the example.  

\subsection{Example 1\label{sect:Jorris2D}}
Consider the following chance constrained variation of a deterministic optimal control problem from Ref~\cite{Jorris2}.  Minimize the cost functional
\begin{equation}\label{eq:costvers2}
J = t_f,
\end{equation} 
subject to the dynamic constraints
\begin{equation} \label{eq:dynvers2}
\begin{array}{ccc}
\dot x (t) & = & V \cos \theta (t), \\ 
\dot y (t) & = & V \sin \theta (t), \\ 
\dot \theta (t) & = & \frac{\tan(\sigma_{\max})}{V}u, \\ 
\dot V & = & a,
\end{array}
\end{equation}
the boundary conditions
\begin{equation}\label{eq:boundvers2}
  \begin{array}{cccccc}
 \ x(0) & = & -1.385, & x(t_f) & = & 0.516, \\ 
 \ y(0) & = & 0.499, & y(t_f) & = & 0.589, \\ 
 \theta (0) & = & 0, & \theta(t_f) & = & \textrm{free}, \\
 V (0) & = & 0.293, & V(t_f) & = & \textrm{free}, \\ 
\end{array}
 \end{equation}
 the control bounds
 \begin{equation}\label{eq:contvers2}
 -1 \leq u \leq 1,
 \end{equation}
 the event constraints
 \begin{equation}\label{eq:eventvers2}
 \big( x(t_i)-x_i,y(t_i)-y_i \big)=\big( 0,0 \big),\quad (i=1,2),
 \end{equation}
and the chance path inequality constraint (keep-out zone constraint)
\begin{equation}\label{eq:CC1}
 P \left( R^2 - \Delta x_{\xi_1}^2-\Delta y_{\xi_2}^2  > \delta \right) \leq \epsilon_d, 
\end{equation}
where $(\Delta x_{\xi_1}, \Delta y_{\xi_2})$ are defined as
\begin{equation}\label{eq:deltas}
\big( \Delta x_{\xi_1},\Delta y_{\xi_2} \big)  = \big( x+\xi_1 -x_c, y+\xi_2-y_c \big).
\end{equation}
The random variables $\xi_1$ and $\xi_2$ have normal distributions of $N(\mu_1,\sigma_1^2)$ and $N(\mu_2,\sigma_2^2)$, respectively.  The parameters for the example problem are provided in Table~\ref{table3}.

\begin{table}[ht]
\caption{Parameters for Example 1.\label{table3}}
\renewcommand{\baselinestretch}{1}\small\normalfont
\centering
\begin{tabular}{| c | c |}
  \hline
  Parameter & Value \\ \hline\hline 
  $(x_c,y_c)$ & $(0.193, 0.395)$\\ \hline
  $R$  & $0.243$ \\ \hline
  $\epsilon_d$ & $0.010$ \\ \hline
  $\delta$ & $0.020$  \\ \hline  
  $a $ & $-0.010$ \\ \hline
  $(x_1, y_1)$ & $(-0.737,0.911)$ \\ \hline
  $(x_2, y_2)$ & $(-0.340,0.297)$ \\ \hline
  $\sigma_{\textrm{max}}$ & $0.349$ \\ \hline
  $ (\mu_1,\mu_2)$ & $(0,0)$ \\ \hline
  $(\sigma_1,\sigma_2)$ & $(0.001,0.0005)$ \\ \hline
\end{tabular}
\end{table}

\subsection{Initialization for Example 1}\label{sect:guesseschoic}

For the initialization of Example 1, the problem is divided into three phases.  The first phase begins at $ ( x(0),y(0)) $ and ends at $(x_1,y_1)$.  The second phase ends at $(x_2,y_2)$, while the third phase ends at the $(x(t_f),y(t_f))$.  The constraints of Eqs.~\eqref{eq:dynvers2}--\eqref{eq:contvers2} and Eq.~\eqref{eq:CC1} are included in every phase.  

\begin{figure}[ht!]
  \centering
  \vspace*{0.25in}
{\includegraphics[height = 2.1in]{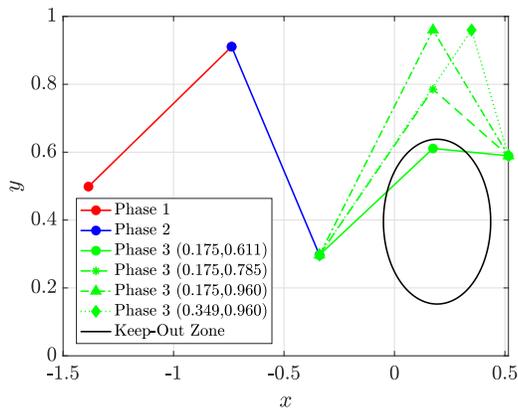}}
\caption{State initial guesses for Example 1.\label{fig:InitialG}}
\end{figure}

As Example 1 is divided into phases, an initial guess of the states and control must be provided for each phase.  Because the deterministic path constraint is active in the third phase and this constraint depends only on $(x,y)$, the initial guess of $(x,y)$ for the third phase affects whether or not the NLP solver converges to solution.  Consequently, as shown in Fig.~\ref{fig:InitialG}, the different initial guesses of $(x,y)$ for the third phase are two line segments connected at the points $(0.175,0.611)$, $(0.175,0.785)$, $(0.175,0.960)$ and $(0.349,0.960)$.  These initial guesses will be referred to, respectively, as initial guesses I, II, III, and IV as shown in Table~\ref{tableInitialGuessesPhase3}.  For all other states, a straight line approximation between the known initial and terminal conditions per phase was applied.  If endpoint conditions were not available, the same constant initial guess that did not violate the constraint bounds for each phase was used.  The control was set as a constant of zero for all three phases.

\begin{table}[ht!]
  \centering
  \caption{Initial guesses used for phase 3 of Example 1 and corresponding to Fig.~\ref{fig:InitialG}.\label{tableInitialGuessesPhase3}}
  \renewcommand{\baselinestretch}{1}\small\normalfont
  \begin{tabular}{|c|c|}\hline
    Initial Guess & Label \\ \hline
    $(0.175,0.611)$ & I \\ \hline
    $(0.175,0.785)$ & II \\ \hline
    $(0.175,0.960)$ & III \\ \hline
    $(0.349,0.960)$ & IV \\  \hline
  \end{tabular}
\end{table}

The chance path constraint is transformed to a deterministic path constraint using the approach of Section~\ref{sect:review}, and evaluating this deterministic path constraint at each collocation point can be computationally intractable when the number of samples is sufficiently large.  In order to ensure computational tractability, the chance path constraint of Eq.~\eqref{eq:CC1} is reformulated as follows~\cite{Keil2}:
\begin{equation}\label{eq:CC2}
\epsilon_d \geq 
\begin{cases}
0, \ & \textrm{if} \ \Delta x^2 + \Delta y^2 \geq (R+b)^2, \\
P \left( R^2 - \Delta x_{\xi_1}^2-\Delta y_{\xi_2}^2  - \delta > 0 \right), &  \textrm{if} \ \Delta x^2 + \Delta y^2 < (R+b)^2,
 \end{cases}
\end{equation}
where $(\Delta x, \Delta y)$ is defined as
\begin{equation}\label{eq:Deltas_noxi}
\big( \Delta x,\Delta y) = \big( x-x_c,y-y_c), \\
\end{equation}
and $b$ is a user defined parameter for determining when the chance path constraint will be evaluated with or without samples.  Due to the size of $R$, $b$ is set equal to $0.05$.  The chance path constraint of Eq.~\eqref{eq:CC1} evaluates to a small number if the distance between $(x,y)$ and $(x_c,y_c)$ is large enough, and this small number indicates that the chance path constraint is inactive.  Consequently, if the distance between $(x,y)$ and $(x_c,y_c)$ is larger than $R+b$, the chance path constraint is taken to be inactive and will be set equal to an arbitrary constant less than $\epsilon_d$ (which in this case is zero).   Conversely, if the distance between $(x,y)$ and $(x_c,y_c)$ is smaller than $R+b$, the chance path constraint will be evaluated using samples.  Therefore, when the chance path constraint is transformed to a deterministic path constraint, this deterministic path constraint will only be evaluated using samples at a subset of the collocation points, thus improving  computational tractability.  

\subsection{Setup for Optimal Control Software $\mathbb{GPOPS-II}$ \label{sect:setup}}

Example 1 was solved using the $hp$-adaptive Gaussian quadrature collocation \cite{Garg1,Garg2,Garg3,Patterson2015,Liu2015,Liu2018,Darby2,Darby3,Francolin2014a} optimal control software $\mathbb{GPOPS-II}$~\cite{Patterson2014} together with the NLP solver SNOPT~\cite{Gill1,Gill2} (using a maximum of $500$ NLP solver iterations) and the mesh refinement method of Ref.~\cite{Liu2018}.  All derivatives required by the NLP solver were obtained using sparse central finite-differencing~\cite{Patterson2012}.  The initial mesh consisted of $10$ mesh intervals with four collocation points each.  Next, constant bandwidths for each kernel were determined using the MATLAB$^{\textrm{\textregistered}}$ function \textsf{ksdensity}~\cite{bowman1,silverman1}.  Furthermore, the method of Neal~\cite{MCMCMethods,Neal2} was used to obtain $50,000$ MCMC samples per run.  Finally, twenty runs for each kernel were performed using a 2.9 GHz Intel$^{\textrm{\textregistered}}$ Core i9 Macbook Pro running Mac OS-X version 10.13.6 (High Sierra) with 32 GB 2400 MHz DDR4 RAM using MATLAB$^{\textrm{\textregistered}}$ version R2018a (build 9.4.0.813654) 

Because a new set of samples is generated for each run, it is not guaranteed that the NLP solver will converge to a solution on every run.  Therefore, consecutive runs are performed to determine the reproducibility of the results.  Consequently, the same number of mesh refinements must be applied per run.  In order to ensure consistent results, the number mesh refinement iterations is limited to two.  It is noted that, from trial and error, twenty runs was found to be sufficient to determine if there were issues of reproducibility, such as a run not converging, because these issues would surface for at least one of the twenty runs.

\subsection{Results and Discussion for Example 1: Without a Warm Start}\label{sect:Jorris2Ddiscuss}

In this section, results are provided for solving Example 1 using the approach of Section~\ref{sect:review} without a warm start, and with the following three kernels: the Split-Bernstein~\cite{Keil2} kernel, the Epanechnikov kernel~\cite{Epanech1} with a bias equal to the bandwidth, and the Gaussian kernel with a bias equal to three times the bandwidth.  The Gaussian kernel was chosen despite not satisfying the requirements of a biased kernel, so that solutions could be obtained using a smooth kernel~\cite{Keil2}.

\begin{figure}[ht]
  \centering
  \vspace*{0.25in}
\subfloat[States.]{\includegraphics[height = 2.1in]{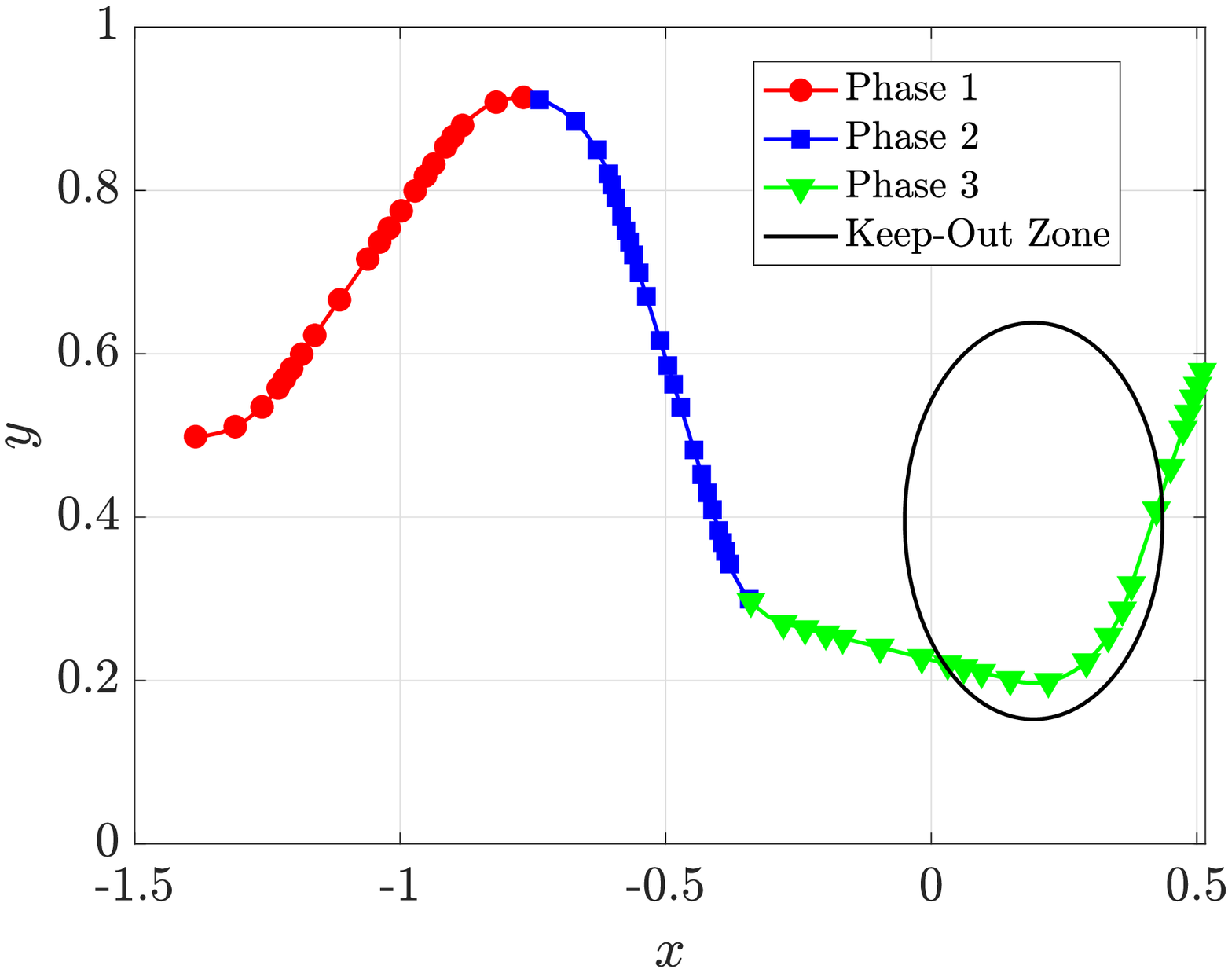}}
~~~~\subfloat[Control.]{\includegraphics[height = 2.1in]{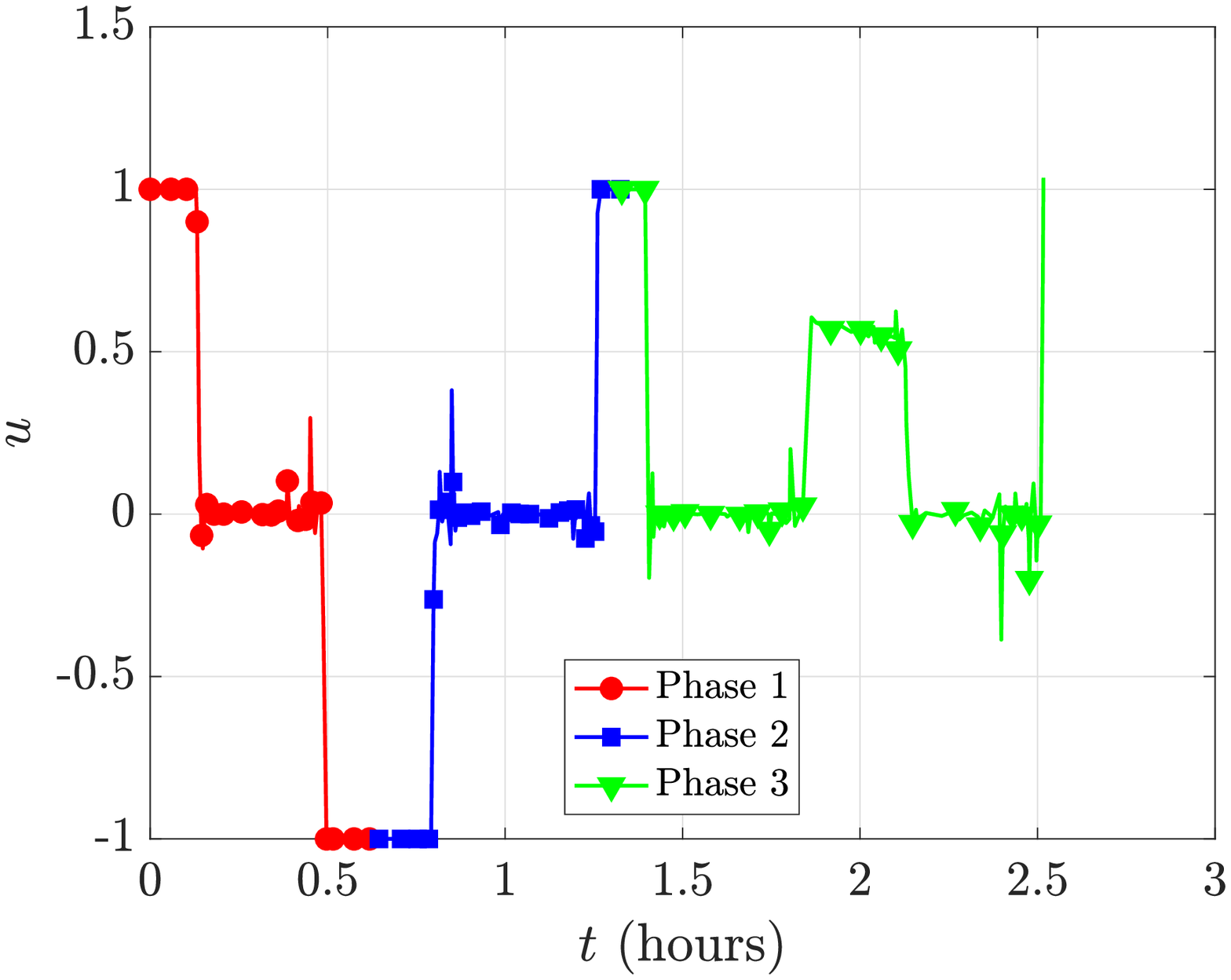}}
\caption{Higher cost solution for Example 1.}
\label{fig:LocalResult}
\end{figure}

\begin{figure}[ht]
  \centering
  \vspace*{0.25in}
\subfloat[States.]{\includegraphics[height = 2.1in]{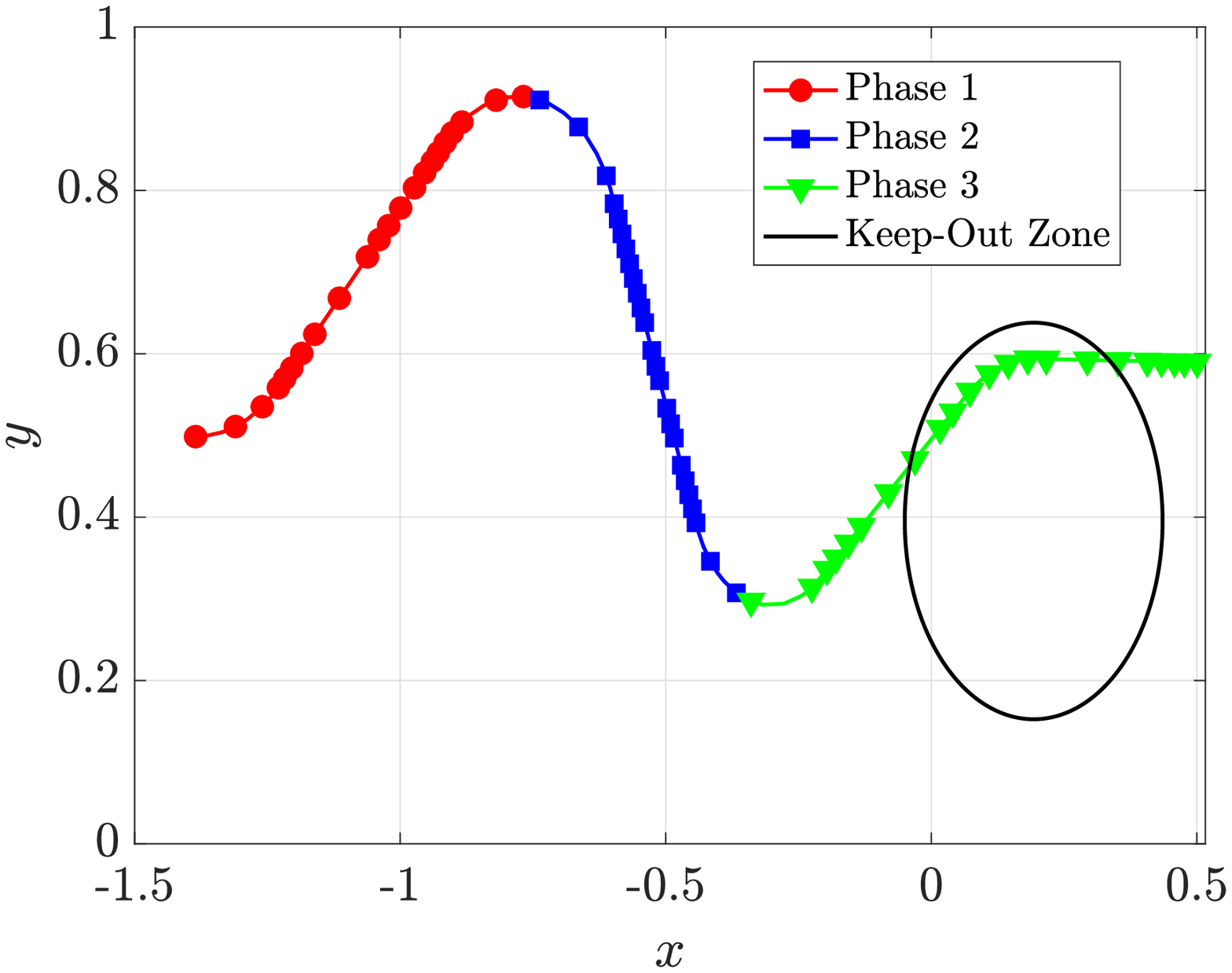}}
~~~~\subfloat[Control.]{\includegraphics[height = 2.1in]{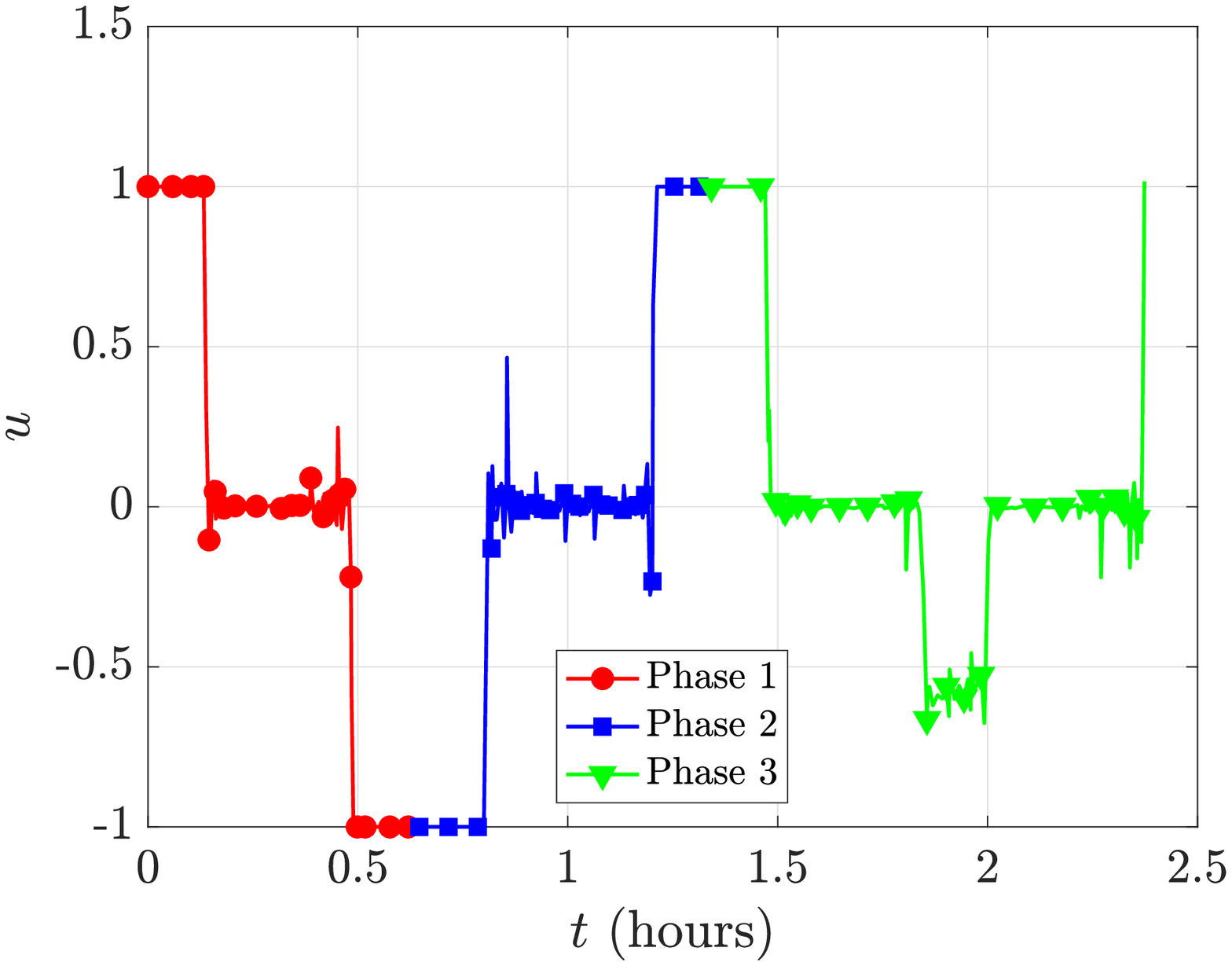}}
\caption{Lower cost solution for Example 1.}
\label{fig:OptResult}
\end{figure}

For Example 1, the NLP solver could either converge to a higher cost or lower cost solution, or not converge.  Figures~\ref{fig:LocalResult} and~\ref{fig:OptResult} show the two different possible solutions obtained using the Split-Bernstein kernel.   It is noted that, for the Gaussian kernel, an infeasible solution was obtained for one run with initial guess II.  This was the only infeasible solution obtained for all of the runs.  Additionally, for Example 1, Tables~\ref{tableSB_naive}--\ref{tableG_naive} contain the results of twenty runs using, respectively, the Split-Bernstein, Epanechnikov, and Gaussian kernels.  For Tables~\ref{tableSB_naive}--\ref{tableG_naive}, $C$ is the number of times the NLP solver converged and $H$ is the number of times the NLP solver converged to the higher cost solution.  Additionally, $\mu_T$, $T_{\min}$, and $T_{\max}$ are, respectively, the average, minimum, and maximum of the execution times for $\mathbb{GPOPS-II}$ obtained for all the runs.  Comparing the results shown in Tables~\ref{tableSB_naive}--\ref{tableG_naive}, it is seen that a large percentage of the runs either converge to the higher cost solution or do not converge.  The best convergence results were for the Epanechnikov kernel with the initial guess $(0.349,0.960)$.  Furthermore, the run times for all three kernels are high regardless of the initial guess, where the best run times were for the Gaussian kernel with the initial guess $(0.175,0.960)$.  Thus, both the kernel and initial guess have an impact on convergence of the NLP solver and run time.  

\begin{table}[ht]
  \centering
  \caption{Results for Example 1 Using the Split Bernstein, Epanechnikov, and Gaussian Kernels.\label{tableExample1Results}}
  \renewcommand{\baselinestretch}{1}\small\normalfont
  \subfloat[Results for Split Bernstein Kernel.\label{tableSB_naive}]{
    \begin{tabular}{| c || c | c | c | c |}
      \hline
      \backslashbox{Quantity}{Initial Guess} & I & II & III & IV  \\ \hline 
      $C$ & $17$ & $14$ & $18$ & $18$ \\ \hline
      $H$ & $7$ & $7$ & $19$ & $3$ \\ \hline
      $\mu_T$ (s) & $314.6$  & $165.1$  & $129.1$  & $ 192.4$ \\ \hline
      $T_{\min}$ (s) & $52.85$  & $43.05$  & $47.69$  & $59.61$  \\ \hline 
      $T_{\max}$ (s) & $2465.8$  & $398.3$  & $1090.4$  & $1916.8$  \\ \hline
    \end{tabular}
  }

  \subfloat[Results for Epanechnikov Kernel.\label{tableE_naive}]{
  \begin{tabular}{| c || c | c | c | c |}
    \hline
    \backslashbox{Quantity}{Initial Guess} & I & II & III & IV  \\ \hline 
    $C$ & $18$ & $20$ & $20$ & $19$ \\ \hline
    $H$ & $3$ & $16$ & $19$ & $4$ \\ \hline
    $\mu_T$ (s) & $346.0$  & $363.5$  & $159.2$  & $261.1$ \\ \hline
    $T_{\min}$ (s) & $84.64$  & $80.88$  & $51.94$  & $64.43$  \\ \hline 
    $T_{\max}$ (s) & $2656.2$  & $1375.0$  & $517.6$  & $1954.8$  \\ \hline
  \end{tabular}
}

\subfloat[Results for Gaussian Kernel.\label{tableG_naive}]{
  \begin{tabular}{| c || c | c | c | c |}
    \hline
    \backslashbox{Quantity}{Initial Guess} & I & II & III & IV  \\ \hline 
    $C$ & $17$ & $19$ & $20$ & $19$ \\ \hline
    $H$ & $3$ & $3$ & $20$ & $9$ \\ \hline
    $\mu_T$(s) & $197.4$  & $287.8$  & $45.98$  & $98.23$ \\ \hline
    $T_{\min}$ (s) & $55.38$  & $41.38$  & $34.61$  & $44.27$  \\ \hline 
    $T_{\max}$ (s) & $980.7$  & $3060.8$  & $85.07$  & $541.7$  \\ \hline
  \end{tabular}
}
\end{table}

This example demonstrates that, without a warm start, solving a complex CCOCP using the approach of Section~\ref{sect:review} can be computationally challenging.  In particular, these computational challenges include the inability to solve the NLP and large computation times.  Moreover, it is noted that these computational issues are affected by the choice of the kernel and initial guess.  In order to overcome these computational issues, in Section \ref{sect:tech} a warm start method is developed for solving CCOCPs.

\section{Warm Start Method\label{sect:tech}}

In this section, a warm start method is developed for efficiently solving CCOCPs using the approach of Section~\ref{sect:review}.  The warm start method consists of three components that are designed to aid the NLP solver in converging to a solution.  Once the NLP solver has converged, the components are applied to efficiently cycle through mesh refinement iterations.  The three components are: (1) bandwidth tuning (Section \ref{sect:bandwidth-tuning}); (2) kernel switching (Section \ref{sect:kernel-switching}); and (3) sample size increase (Section \ref{sect:sample-size-increasing}).  Section~\ref{sect:sumpre} summarizes the warm start method.

\subsection{Component 1:  Tuning the Bandwidth\label{sect:bandwidth-tuning}}

The first component of the warm start method is tuning the bandwidth.  The need for tuning the bandwidth arises from the deterministic constraint obtained by transforming the chance constraint using the approach of Section~\ref{sect:review} being difficult for the NLP solver to evaluate.  Increasing the size of the bandwidth will improve the starting point for the NLP solver, and the NLP solver is more likely to converge with this better starting point.  The solution obtained using the larger bandwidth will, however, have a higher cost than the solution obtained using the original bandwidth.  It is noted that using a larger starting bandwidth that is later reduced to the original bandwidth increases the likelihood that the NLP solver will converge.  Moreover, the solution ultimately obtained will be for the original non-smooth constraint.  

In this paper, the following approach is used to tune the bandwidth.  First, the starting bandwidth is the original bandwidth multiplied by a constant $w \geq 1$, where the original bandwidth is obtained using the MATLAB function \textsf{ksdensity}.  When the mesh refinement error is less that a user chosen parameter $\phi$, $w$ is set to unity.  For tuning the bandwidth, $w$ is started at unity and increased over a series of trial runs of $\mathbb{GPOPS-II}$ until either the NLP solver converges on every run, or the NLP solver can no longer converge to the same solution as when the original bandwidth was used.  It is noted that this approach for tuning the bandwidth is different from that of Ref.~\cite{Keil2} in that Ref.~\cite{Keil2} does not choose the starting bandwidth relative to the original bandwidth.  Finally, it is noted that this first component of the warm start method requires the starting bandwidth to be tuned separately for each kernel.

\subsection{Component 2:  Kernel Switching\label{sect:kernel-switching}}

The second component of the warm start method is switching the kernel.  The purpose of kernel switch is that, even with bandwidth tuning, the likelihood that the NLP solver will not converge is higher if certain kernels are chosen.  Conversely, via the choice of an appropriate kernel, tuning the bandwidth improves the chances that the NLP solver will converge.  Thus, starting with an appropriate combination of bandwidth and kernel and later switching to the desired kernel improves the chances that the NLP solver will converge, even if the second kernel would have resulted in divergence of the NLP solver if it had been used at the outset.  It is noted that an additional benefit to this switch is that kernels, like the Gaussian kernel, that do not satisfy the criteria for a biased KDE from Section~\ref{sect:review}, can still be applied as the starting kernel.  Only the desired kernel must satisfy this criteria.  

In the method of this paper, the kernel switch is performed as follows.  First, a starting bandwidth and kernel are chosen by trial runs.  Next, the kernel is switched when the mesh refinement error is below $\phi$.  Thus, the bandwidth and kernel are updated simultaneously.

\subsection{Component 3:  Incrementally Increasing Sample Size\label{sect:sample-size-increasing}}

The third component of the warm start method is an approach for incrementally increasing the sample set size.  The purpose of incrementally increasing the sample set size is that it is computationally expensive to evaluate the deterministic constraint obtained by transforming the chance constraint using the approach of Section~\ref{sect:review} when the sample size is large.  If, on the other hand, the number of samples is reduced, the computational effort required by the NLP solver is also reduced~\cite{Roycet1}.  Note, however, that by reducing the sample size, it is no longer possible to satisfy the bound on the deterministic constraint as described in Section~\ref{sect:review}.  Conversely, if the number of samples is increased incrementally from a small amount to the total number of samples through a series of mesh refinement iterations, the computational expense is reduced while ultimately satisfying the bound on the deterministic constraint.  

In this paper, the following approach is used to incrementally increase the sample set size.  First, a small number of samples is selected as the starting sample set.  Next, when the mesh refinement error drops below the user-specified value $\phi$, the number of samples is increased by a user-specified amount.  Thus, the first increase in sample size occurs when the bandwidth and kernel are updated.  After the first increase in the number of samples, the sample size is incrementally increased on every subsequent mesh refinement iteration until the full sample size is reached.  It is noted that, because the increments for increasing the samples are tied to the mesh refinement iterations, using an excessive number of increments may lead to a need for extra mesh refinement iterations.  Moreover, for every set of samples, a different bandwidth must be generated using $\textsf{ksdensity}$.  Consequently, the starting bandwidth will be the bandwidth for the smallest sample set multiplied by $w$ as described in Section \ref{sect:bandwidth-tuning}.  

\subsection{Summary of Warm Start Method}\label{sect:sumpre}

The three components of the warm start method are changing the bandwidth, switching the kernel, and an approach for incrementally increasing the sample set size.  These components increase the chances of the NLP solver converging while reducing run time, regardless of the choice of kernel.  Also, the sensitivity to the initial guess will be reduced by applying an appropriate starting bandwidth, kernel, and subset of samples.  The three components are combined into the warm start method that is presented below.
 
\begin{shadedframe}
\vspace{-10pt}
\begin{center}
 \shadowbox{\bf Warm Start Method for Solving CCOCPs}
\end{center}
\begin{enumerate}[{\bf Step 1:}]
\item Determine bandwidths for $2000$ and $10,000$ samples from full sample set, as well as for the full sample set.  
\item Choose a constant $w$ and kernel pair.\label{step:multiplier}
\begin{enumerate}[{\bf (a):}]
\item Choose a trial constant $w_i$ and kernel.
\item Perform 10-20 runs for up to two mesh refinement iterations, using $2000$ samples from full sample set.\label{step:refine}
\item If NLP solver converges on all runs, set $w = w_i$.  Otherwise, choose $w_{i+1} > w_{i}$ and possibly change the kernel.  Return to {\bf~\ref{step:refine}}.
\end{enumerate}
\item Run problem through optimal control software for a series of mesh refinement iterations with $2000$ samples from full sample set. \label{step:runfull}
\item For the first mesh refinement iteration when mesh error decreases to $\phi$: set $w = 1$, change to $10,000$ samples from full sample set, update the bandwidth, and switch the kernel.
\item On the following mesh refinement iterations, change the sample set to the full set of samples and update the bandwidth.
\end{enumerate}
\end{shadedframe}

\section{Solution to Example 1 Using Warm Start Method \label{sect:discusstech}}

Example 1 is now re-solved using biased KDEs and LGR collocation together with the warm start method of Section \ref{sect:tech}.  For Example 1, the values $\phi = 5 \times 10^{-5}$ and $w = 100$ are used, and the Split-Bernstein kernel is the starting kernel.  In Section~\ref{sect:lowmesh}, Example 1 is solved using three mesh refinement iterations such that the final mesh refinement is performed using the full sample set, in order to provide a fair comparison with the results obtained without a warm start as given in Section \ref{sect:Jorris2Ddiscuss} (where it is noted that two mesh refinement iterations were used to obtain the results shown in Section \ref{sect:Jorris2Ddiscuss}).  In Section~\ref{sect:highmesh}, Example 1 is solved with enough mesh refinement iterations to reach mesh convergence, along with a deterministic version of Example 1.

\subsection{Limited Number of Mesh Refinement Iterations\label{sect:lowmesh}}

Recall from Section \ref{sect:Jorris2Dnaive} that a maximum of two mesh refinement iterations was allowed when solving Example 1 without a warm start.  To demonstrate the improvement of using a warm start while providing a fair comparison with the results obtained in Section \ref{sect:Jorris2Dnaive}, in this section, Example 1 is solved using the approach of Section~\ref{sect:review}, with the warm start method of Section~\ref{sect:tech} and a maximum of three mesh refinement iterations.  Tables~\ref{tableSB_threeiter}--\ref{tableG_threeiter} show the results obtained using the warm start method for, respectively, the Split-Bernstein, Epanechnikov and Gaussian kernels.  The results show that, with the warm start method, the NLP solver converges to the lower cost solution for all runs, as compared to converging to the higher cost solution or not converging without a warm start, as shown previously in Tables~\ref{tableSB_naive}--\ref{tableG_naive}.  Also, the results obtained using the warm start method indicate that convergence of the NLP solver was not affected by the kernel or initial guesses.

Furthermore, the run times using all three kernels are much lower when a warm start is included, when compared with the results of Section~\ref{sect:Jorris2Dnaive}.  In addition, the run times are similar regardless of the choice of the kernel or the initial guess.  This last observation indicates that the computation time is not affected significantly by the choice of the kernel or the initial guess.  To see the differences between the results with and without the warm start method, Table~\ref{perform} provides the percentage increase in computational performance when solving Example 1 with a warm start relative to not including a warm start.  The results in Table~\ref{perform} show that the most significant difference between including and not including the warm start method occurs in the maximum and average computation times.  Additionally, the difference in the performance increase between the three kernels is insignificant.

\begin{table}[ht]
  \centering
  \caption{Results for Example 1 with a warm start. \label{tableExample1ResultsWarmStart}}
  \renewcommand{\baselinestretch}{1}\small\normalfont
  \subfloat[Results for Split Bernstein kernel.\label{tableSB_threeiter}]{
    \begin{tabular}{| c || c | c | c | c |}
      \hline
      \backslashbox{Quantity}{Initial Guess} & I & II & III & IV  \\ \hline 
      $C$ & $20$ & $20$ & $20$ & $20$ \\ \hline
      $H$ & $0$ & $0$ & $0$ & $0$ \\ \hline
      $\mu_T$ (s) & $29.79$  & $32.15$  & $37.23$  & $30.54$ \\ \hline
      $T_{\min}$ (s) & $23.30$  & $23.62$  & $22.82$  & $18.68$  \\ \hline 
      $T_{\max}$ (s) & $49.92$  & $68.59$  & $83.15$  & $85.10$  \\ \hline
    \end{tabular}
  }

  \subfloat[Results for Epanechnikov kernel.\label{tableE_threeiter}]{
    \begin{tabular}{| c || c | c | c | c |}
      \hline
      \backslashbox{Quantity}{Initial Guess} & I & II & III & IV  \\ \hline 
      $C$ & $20$ & $20$ & $20$ & $20$ \\ \hline
      $H$ & $0$ & $0$ & $0$ & $0$ \\ \hline
      $\mu_T$ (s) & $31.71$  & $38.34$  & $33.27$  & $33.14$ \\ \hline
      $T_{\min}$ (s) & $18.70$  & $28.41$  & $27.36$  & $27.59$  \\ \hline 
      $T_{\max}$ (s) & $51.80$  & $71.13$  & $53.85$  & $49.42$  \\ \hline
    \end{tabular}
  }

\subfloat[Results for Gaussian kernel.\label{tableG_threeiter}]{
  \begin{tabular}{| c || c | c | c | c |}
    \hline
    \backslashbox{Quantity}{Initial Guess} & I & II & III & IV  \\ \hline 
    $C$ & $20$ & $20$ & $20$ & $20$ \\ \hline
    $H$ & $0$ & $0$ & $0$ & $0$ \\ \hline
    $\mu_T$ (s) & $25.67$  & $21.43$  & $21.32$  & $23.29$ \\ \hline
    $T_{\min}$ (s) & $18.67$  & $15.68$  & $15.43$  & $17.72$  \\ \hline 
    $T_{\max}$ (s) & $75.86$  & $27.94$  & $26.22$  & $30.96$  \\ \hline
  \end{tabular}
  }
\end{table}

\begin{table}[h!]
\small
\caption{Increase in performance using the warm start method.\label{perform}}
\renewcommand{\baselinestretch}{1}\small\normalfont
\centering
\begin{tabular}{| c || c | c | c | c |}
  \hline
  Initial Guess & I & II & III & IV \\ \hline
  Split-Bernstein $\mu_T$ & $90.53 \% $ &$80.53 \%$ & $ 71.15 \%$ & $84.12 \%$  \\ \hline
  Split-Bernstein $T_{\min}$ & $ 55.91 \% $ &$45.13 \%$ & $52.14 \%$ & $68.66 \%$  \\ \hline
  Split-Bernstein $T_{\max}$ & $97.96 \% $ &$82.78 \%$ & $92.37 \%$ & $95.56 \%$   \\ \hline
  Epanechnikov $\mu_T$ & $90.84 \% $ &$89.45 \%$ & $79.11 \%$ & $87.31 \%$ \\ \hline
  Epanechnikov $T_{\min}$ &$77.91 \%$ &$64.87 \%$ &$47.32 \%$ & $57.18 \%$ \\ \hline
  Epanechnikov $T_{\max}$ & $98.05 \%$ & $94.83 \%$ & $89.60 \%$ & $97.47 \%$ \\ \hline
  Gaussian $\mu_T$ & $87.00 \% $ &$92.55 \%$ & $53.63 \%$ & $76.29 \%$  \\ \hline
  Gaussian $T_{\min}$ & $66.28 \% $ &$62.10 \%$ & $55.41 \%$ & $55.97 \%$  \\ \hline
  Gaussian $T_{\max}$ & $92.26 \% $ &$99.09 \%$ & $69.17 \%$ & $94.28 \%$  \\ \hline
\end{tabular}
\end{table}

\clearpage

\subsection{Unlimited Number of Mesh Refinement Iterations\label{sect:highmesh}}

To further demonstrate the effectiveness of the warm start method of Section~\ref{sect:tech}, Example 1 is now solved using biased KDEs with the warm start method, but with no limit on the number of mesh refinement iterations to reach a user-specified mesh refinement error tolerance.   The limit on the number of mesh refinement iterations is removed because it was found in Section \ref{sect:lowmesh} that the NLP solver converged on every run when the number of mesh refinement iterations was limited.   For the analysis of solving Example 1 with the warm start method and unlimited mesh refinement iterations, a deterministic formulation of Example 1 is also presented, in order to compare the chance constrained solutions to deterministic solutions from the literature.  The deterministic and chance constrained formulations of Example 1 are the same with the exception that the chance constraint of Eq.~\eqref{eq:CC2} is replaced by the following deterministic path constraint
\begin{equation}\label{eq:ex2path}
  R^2 - \Delta x^2-\Delta y^2  \leq 0.
\end{equation}

\begin{figure}[ht!]
  \centering
  \vspace*{0.25in}
\subfloat[Chance Constrained.]{\includegraphics[height = 2.1in]{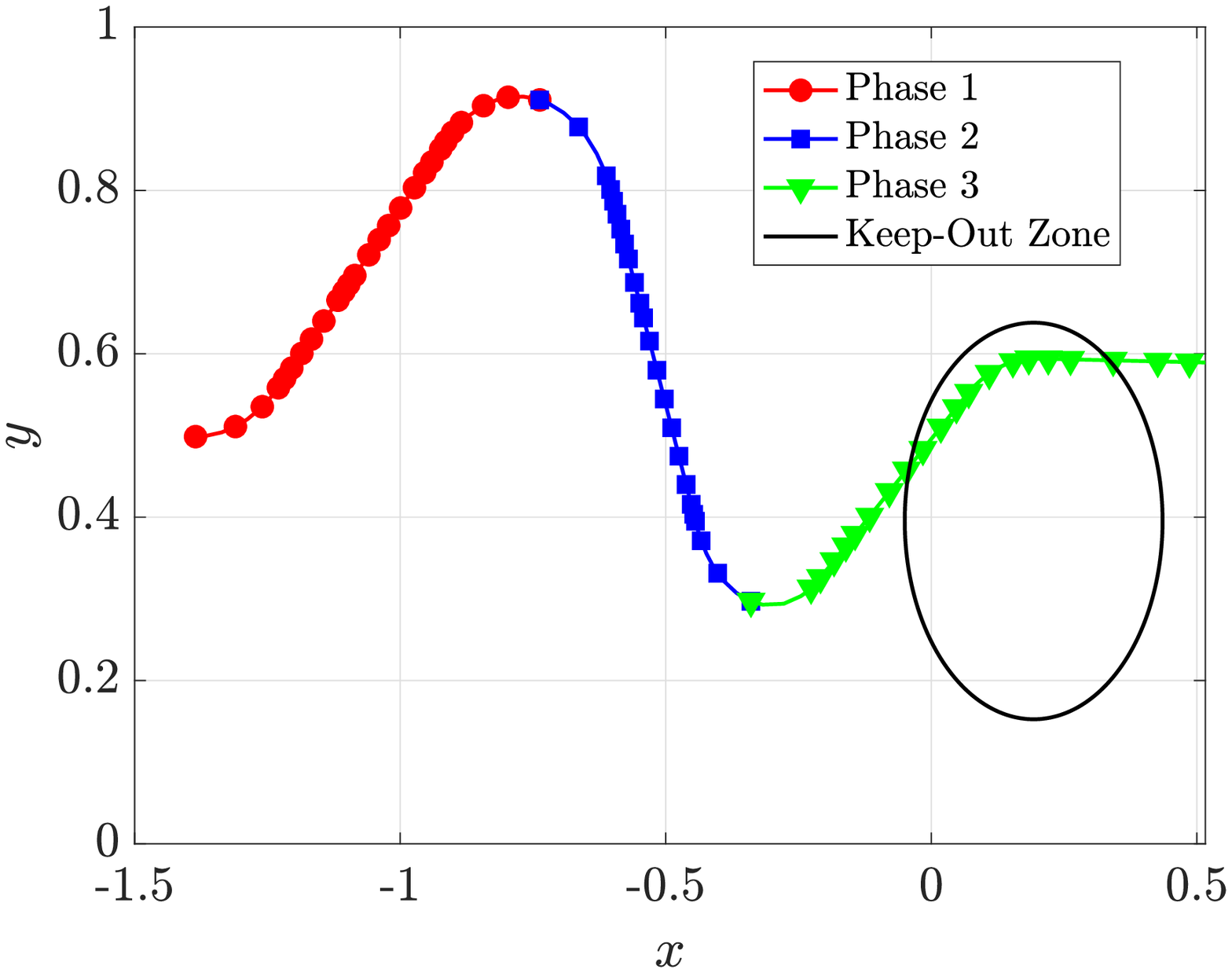}}
~~~~\subfloat[Deterministic.]{\includegraphics[height = 2.1in]{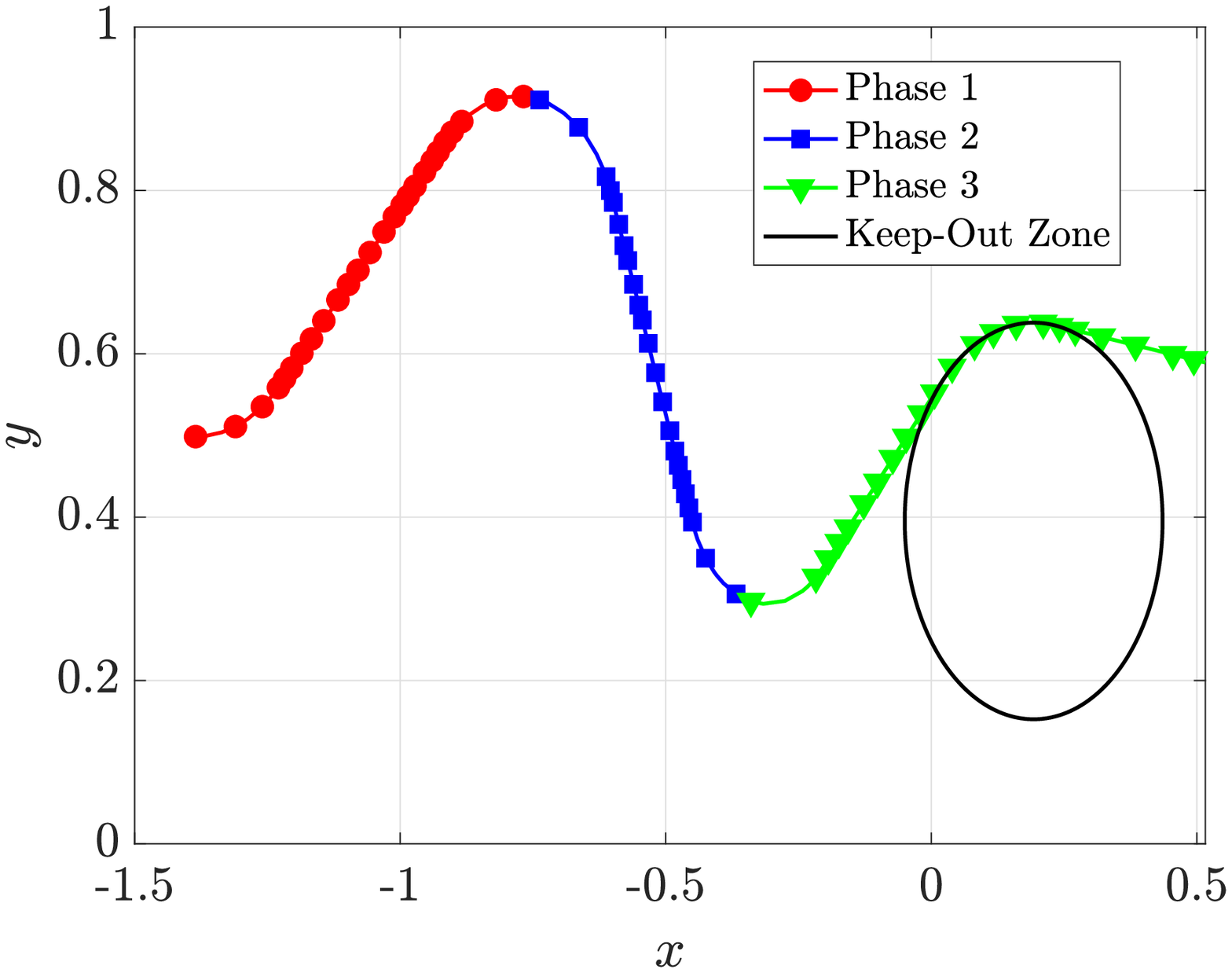}} \\
\subfloat[Chance Constrained.]{\includegraphics[height = 2.1in]{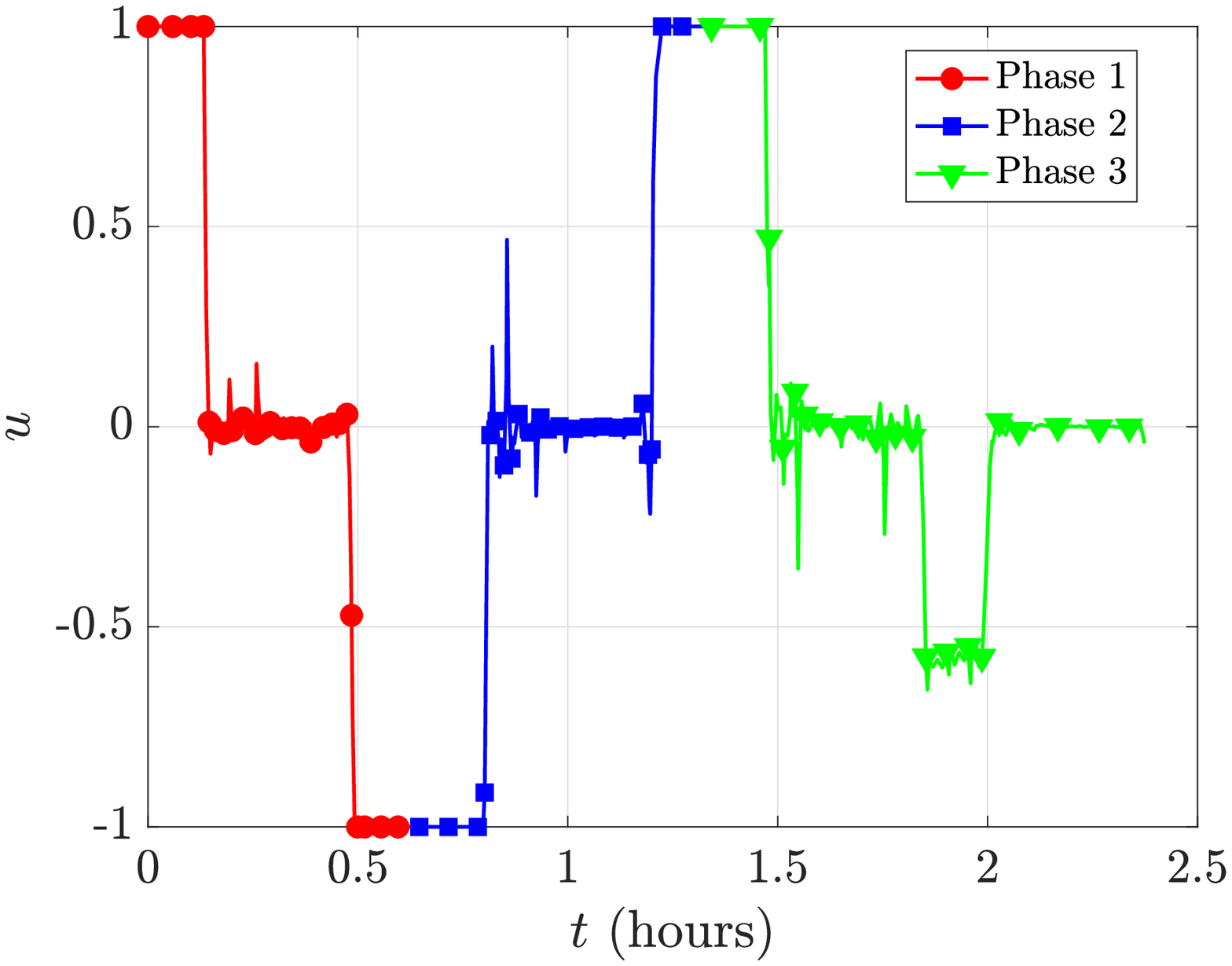}} 
~~~~\subfloat[Deterministic.]{\includegraphics[height = 2.1in]{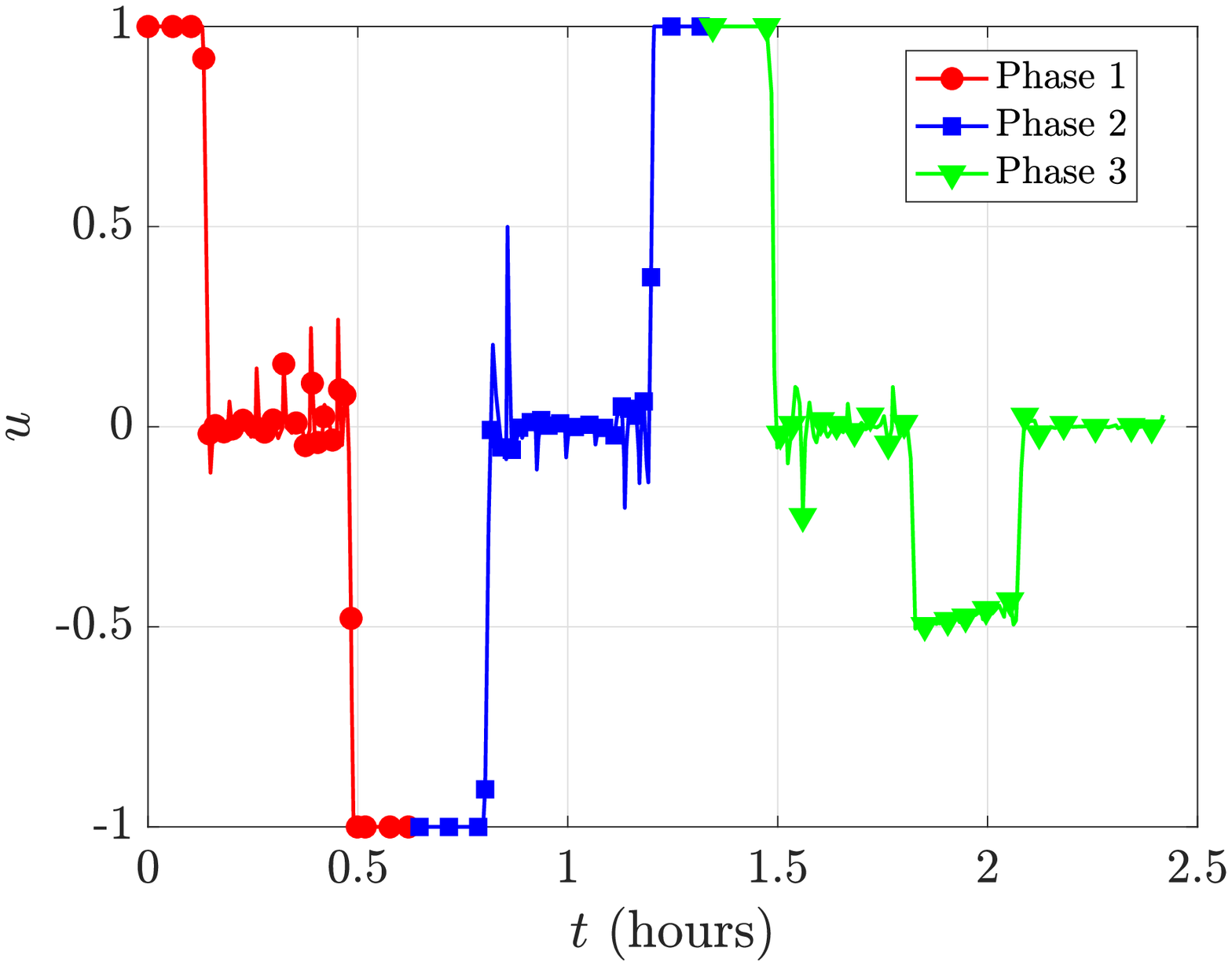}} \\
\caption{Solution for Example 1 and deterministic variation of Example 1.}
\label{fig:Jorris2DEpanech}
\end{figure}

The solution to the chance constrained version of Example 1 using the Epanechnikov kernel is shown in Fig.~\ref{fig:Jorris2DEpanech} alongside the solution to the deterministic formulation of Example 1, where it is seen that the chance constrained and deterministic solutions are similar.  Next, Table~\ref{table_ex1_final} compares the following results: (1) results obtained using the approach of Section~\ref{sect:review} for twenty runs of the chance constrained version of the example [Eqs.~\eqref{eq:costvers2}--\eqref{eq:eventvers2} and~\eqref{eq:CC2}] applying two different kernels, and (2) results obtained for twenty runs of the deterministic formulation of the example [Eqs.~\eqref{eq:costvers2}--\eqref{eq:eventvers2} and \eqref{eq:ex2path}].  For Table~\ref{table_ex2_final}, the quantities $\mu_{J^*}$ and $\sigma_{J^*}$ are the average and standard deviations, respectively, of the optimal cost obtained over all of the runs.  It is noted that only the Split-Bernstein and Epanechnikov kernel were used to obtain results, because the Gaussian kernel was only applied when the number of mesh refinements were restricted in order to determine if the computational challenges were mitigated by using a smooth kernel.  Now that the kernel affects have been reduced, there is no longer a need for a smooth kernel, particularly when the kernel does not satisfy the criteria for a biased KDE from Section~\ref{sect:review}.

The results indicate that the average optimal cost obtained using the chance constrained formulation was lower than the deterministic optimal cost.  The reason for this difference in cost is that the deterministic keep out zone path constraint is designed so that the path of $(x,y)$ can be outside, or on the boundary of keep out zone of radius $R$.  For the chance constrained formulation, a one percent chance of risk violation ($\epsilon_d = 0.01)$ is allowed so that the path of $(x,y)$ can now be a $\delta$ distance radially inside the keep out zone.  As a result, the $(x,y)$ path shown in Fig.~\ref{fig:Jorris2DEpanech} is shorter for the chance constrained formulation than for the deterministic formulation, and subsequently it will take less time to travel this shorter path.  Thus, because the optimal cost is final time, and it takes less time to travel a shorter path, the average optimal cost for the chance constrained formulation will be lower than for the deterministic formulation.  

Additionally, the run times were higher for the chance constrained formulation than for the deterministic formulation.  The deterministic formulation uses a deterministic keep out zone constraint that is not dependent on samples, as opposed to the chance constrained formulation.  Thus, less computational effort is required to solve the deterministic formulation.

\begin{table}[ht]
\caption{Results for Example 1 with warm start method and unlimited mesh refinement iterations. \label{table_ex1_final}}
\renewcommand{\baselinestretch}{1}\small\normalfont
\centering
\begin{tabular}{| c || c | c | c |} \hline
  & Split-Bernstein & Epanechnikov & Deterministic \\ \hline
  $\mu_{J^*}$ (s) & $8545.052$ & $8545.021$ &  $8704.199$ \\ \hline
  $ \sigma_{J^*}$ (s) & $0.0606$ & $0.0517$ & $0$ \\ \hline
  $\mu_T$ (s) & $39.640$ & $43.445$ &  $1.478$ \\ \hline
  $T_{\min}$ (s) & $21.904$ & $20.813$ & $ 1.430$  \\ \hline
  $T_{\max}$ (s) & $87.333$ & $73.209$ & $ 1.591$  \\ \hline
\end{tabular}
\end{table} 

\section{Example 2 Using Warm Start Method\label{sect:examples}}

To further demonstrate the applicability of the warm start method, in this section the warm start method of Section \ref{sect:tech} is applied to a more complex version of Example 1.  Section~\ref{sect:Jorris3Dpresent} provides both a chance constrained and deterministic formulation of this second example.  Section~\ref{sect:Jorris3Dset} describes the initialization.  Finally, Section~\ref{sect:discusdiff} provides the results obtained when solving the example using the warm start method of Section \ref{sect:tech}.  

\subsection{Example 2}\label{sect:Jorris3Dpresent}

Consider the following chance constrained variation of the deterministic optimal control problem from Ref~\cite{Jorris3}.  Minimize the cost functional
\begin{equation}\label{eq:costvers3}
J = t_f,
\end{equation} 
subject to the dynamic constraints
\begin{equation} \label{eq:dynvers3}
\begin{array}{ccc}
\dot x (t) & = & V \cos \theta (t), \\ 
\dot y (t) & = & V \sin \theta (t), \\ 
\dot h (t) & = & V \gamma (t), \\
\dot V(t) & = & - \frac{B V^2 \exp \big(- \beta r_0 h (1+ c_l^2) \big) }{2 E^*}, \\
\dot \gamma (t) & = & BV \exp(- \beta r_0 h) c_l \cos \sigma - \frac{1}{V} + V, \\
\dot \theta (t) & = & BV \exp(- \beta r_0 h) c_l \sin \sigma \\ 
\end{array}
\end{equation}
the boundary conditions
\begin{equation}\label{eq:boundvers3}
  \begin{array}{cccccc}
    \ x(0) & = & -1.385, & x(t_f) & = & 1.147, \\ 
    \ y(0) & = & 0.499, & y(t_f) & = & 0.534, \\ 
    \ h(0) & = & 0.0190, & h(t_f) & = & 0.0038, \\
    \gamma (0) & = & -0.0262, & \gamma (t_f) & = & \textrm{free}, \\
    V (t_0) & = & 0.927, & V(t_f) & = & \textrm{free}, \\ 
    \theta (0) & = & 0.0698, & \theta(t_f) & = & \textrm{free}, \\
  \end{array}
\end{equation}
the control bounds
\begin{equation}\label{eq:contvers3}
  \begin{array}{ccccc}
    - \frac{\pi}{3} & \leq & \sigma & \leq & \frac{\pi}{3}, \\
    0 & \leq & c_l & \leq & 2,
  \end{array}
\end{equation}
the event constraints
\begin{equation}\label{eq:eventvers3}
  (x(t_i)-x_i,y(t_i)-y_i)=(0,0),\quad (i=1,2),
\end{equation}
the path inequality constraints
\begin{equation}\label{eq:ex3path}
  \begin{array}{c}
    R_1^2 - (x-x_{c,1})^2 - (y-y_{c,1})  \leq  0,  \\
    K\exp \bigg( \beta r_0 \frac{h}{2} \bigg) V^3 -1  \leq  0,
  \end{array}
\end{equation}
and the chance path inequality constraint (keep-out zone constraint)
\begin{equation}\label{eq:CCnofly2}
P \left( R_2^2 - \Delta x_{\xi_1,2}^2-\Delta y_{\xi_2,2}^2  > \delta \right) \leq \epsilon_d,
\end{equation}
where $(\Delta x_{\xi_1,2},\Delta y_{\xi_2,2})$ are defined as 
\begin{equation}\label{eq:finDels}
  \big( \Delta x_{\xi_1,2},\Delta y_{\xi_2,2} \big)  = \big( x+\xi_1 -x_{c,2}, y+\xi_2-y_{c,2} \big).
\end{equation}
The random variables $\xi_1$ and $\xi_2$ have normal distributions of $N(\mu_1,\sigma_1^2)$ and $N(\mu_2,\sigma_2^2)$, respectively.  Furthermore, a deterministic version of Example 2 is also solved in order to compare the solutions for chance constrained and deterministic formulations of Example 2.  The deterministic version of Example 2 is identical to that given in Eqs.~\eqref{eq:dynvers3}--\eqref{eq:ex3path}, with the exception that the chance constraint of Eq.~\eqref{eq:ex3path} is replaced with the following deterministic inequality path constraint  
\begin{equation}\label{eq:ex3path_nofly}
  R_2^2 - \Delta x_2^2-\Delta y_2^2  \leq 0. 
\end{equation}
Finally, the parameters for Example 2 are given in Table~\ref{table_example 2}.

\begin{table}[ht]
\caption{Parameters for Example 2.\label{table_example 2}}
\renewcommand{\baselinestretch}{1}\small\normalfont
\centering
\begin{tabular}{| c | c |}  \hline
  Parameter & Value \\ \hline\hline
  $(x_{c,1},y_{c,1})$ & $(0.008,0.389)$ \\ \hline
  $R_1$ & $0.277$\\ \hline
  $(x_{c,2},y_{c,2})$ & $(1.022, 0.943)$  \\ \hline
  $R_2$  & $0.434$ \\ \hline  
  $\epsilon_d$ & $0.010$ \\ \hline
  $\delta$ & $0.020$ \\ \hline
  $K$ & $0.759$ \\ \hline
  $(x_1,y_1)$  & $(-0.466, 0.594)$ \\ \hline 
  $(x_2,y_2)$  &  $(0.728, 0.580)$ \\ \hline
  $B$ & $942.120$ \\ \hline
  $\beta $ & $1.400 \times 10^{-4}$ \\ \hline
  $r_0$ & $6.408 \times 10^6$ \\ \hline
  $E^*$ & $3.240$ \\ \hline
  $ (\mu_1,\mu_2)$ & $(0,0)$ \\ \hline
  $(\sigma_1,\sigma_2)$ & $(0.0007,0.001)$ \\ \hline
\end{tabular}
\end{table}

\subsection{Initialization for Example 2}\label{sect:Jorris3Dset}

Example 2 is implemented as a four-phase problem.  Phase 1 starts at $(x(0),y(0))$ and terminates when the second path constraint of Eq.~\eqref{eq:ex3path} reaches its boundary.  Next, phases 2, 3, and 4 terminate, respectively, at $(x_1,y_1)$, $(x_2,y_2)$, and $(x(t_f),y(t_f))$.  Furthermore, the constraints of Eqs.~\eqref{eq:dynvers3}--\eqref{eq:ex3path} and Eq.~\eqref{eq:CC2_ex2} are included in every phase.  The initial guess for each phase is a straight line approximation between the known initial and terminal conditions for all states.  For any phase where an endpoint was not available, a constant initial guess that did not violate the constraint bounds was used.  The controls were set as straight line approximations between values within the control bounds.  

In order to maintain computational tractability (see Section~\ref{sect:guesseschoic}), the chance constraint of Eq.~\eqref{eq:CCnofly2} is reformulated as follows~\cite{Keil2}:
\begin{equation}\label{eq:CC2_ex2}
  \epsilon_d \geq 
  \begin{cases}
    0, \ \textrm{if} \ \Delta x_2^2 -\Delta y_2^2 \geq (R_2+b)^2, \\
    \begin{aligned}
      P \left( R_2^2 - \Delta x_{\xi_1,2}^2-\Delta y_{\xi_2,2}^2  - \delta > 0 \right), \\ \ \textrm{if} \ \Delta x_2^2 -\Delta y_2^2 < (R+b)^2,
    \end{aligned}
  \end{cases}
\end{equation}
where $(\Delta x_2,\Delta y_2)$ are defined as
\begin{equation}
  (\Delta x_2,\Delta y_2) = (x-x_{c,2},y-y_{c,2}),
\end{equation}
and $b$ is set equal to $0.07$, due to the size of $R$.

\subsection{Results and Discussion for Example 2}\label{sect:discusdiff}

\begin{figure}[ht!]
  \centering
  \vspace*{0.25in}
\subfloat[Chance Constrained]{\includegraphics[height = 2.1in]{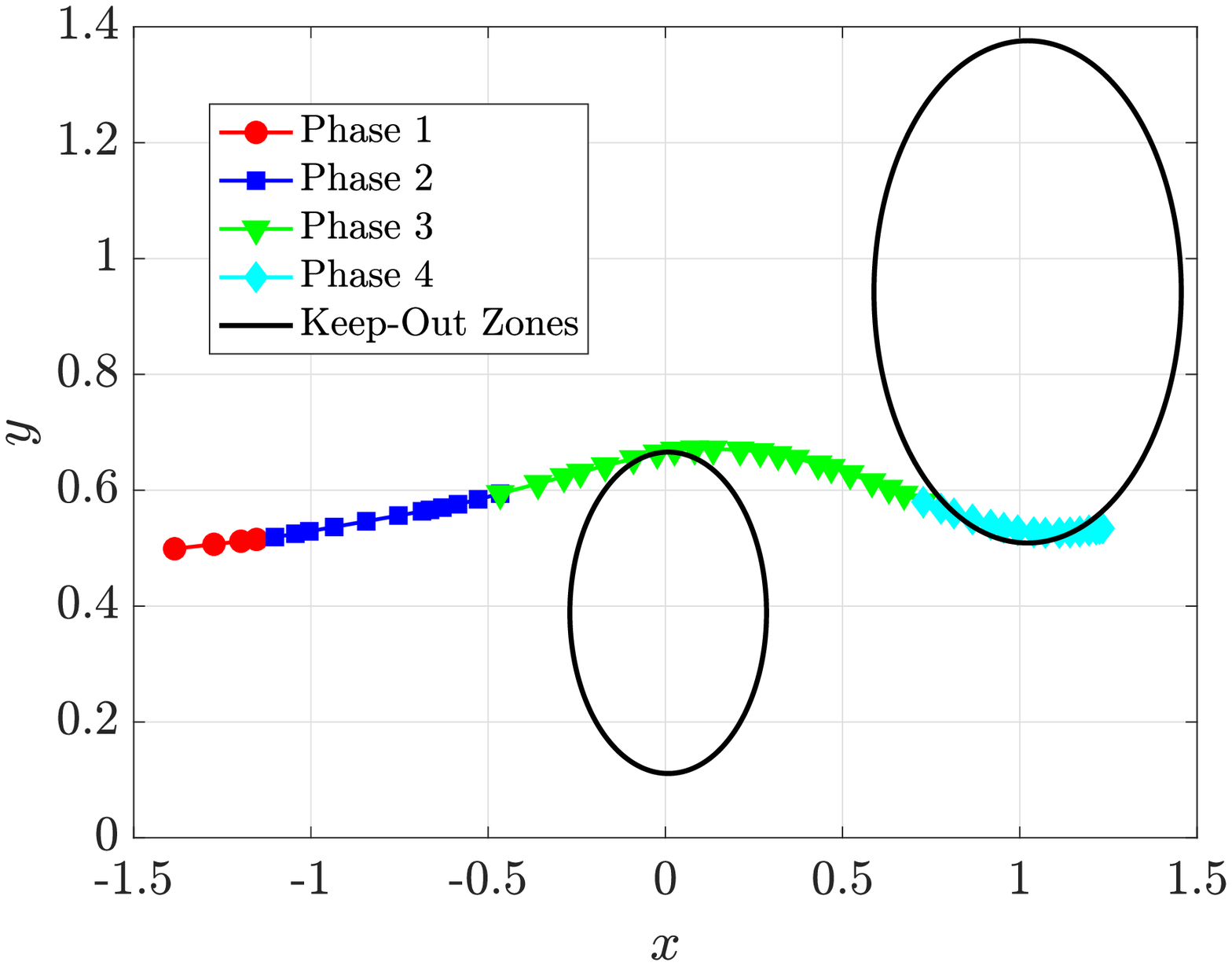}}
~~~~\subfloat[Deterministic]{\includegraphics[height = 2.1in]{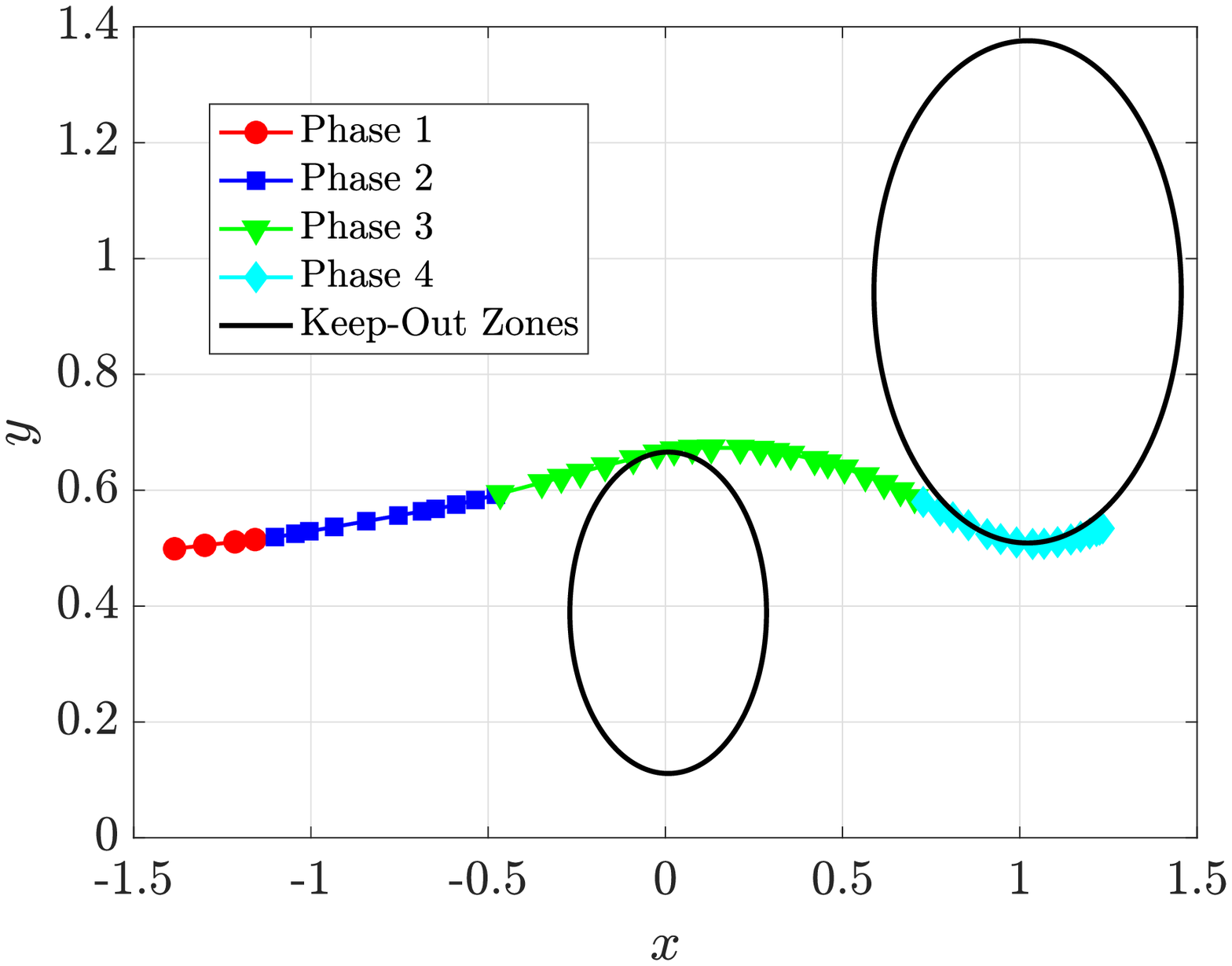}} \\
\caption{States for Example 2 and deterministic variation of Example 2.}
\label{fig:Jorris3DStates}
\end{figure}

\begin{figure}[ht!]
  \centering
  \vspace*{0.25in}
\subfloat[Chance Constrained]{\includegraphics[height = 2.1in]{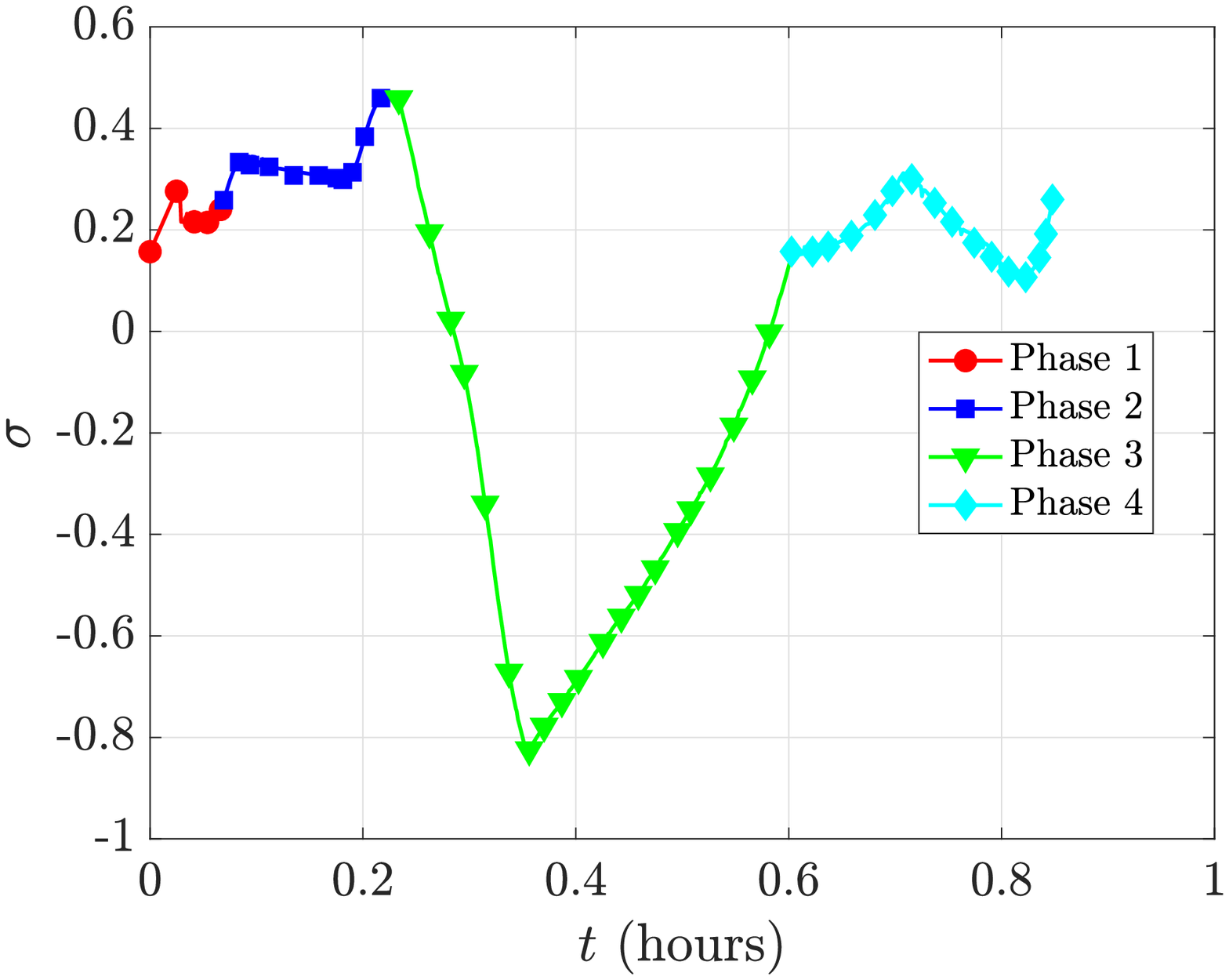}} 
~~~~\subfloat[Deterministic]{\includegraphics[height = 2.1in]{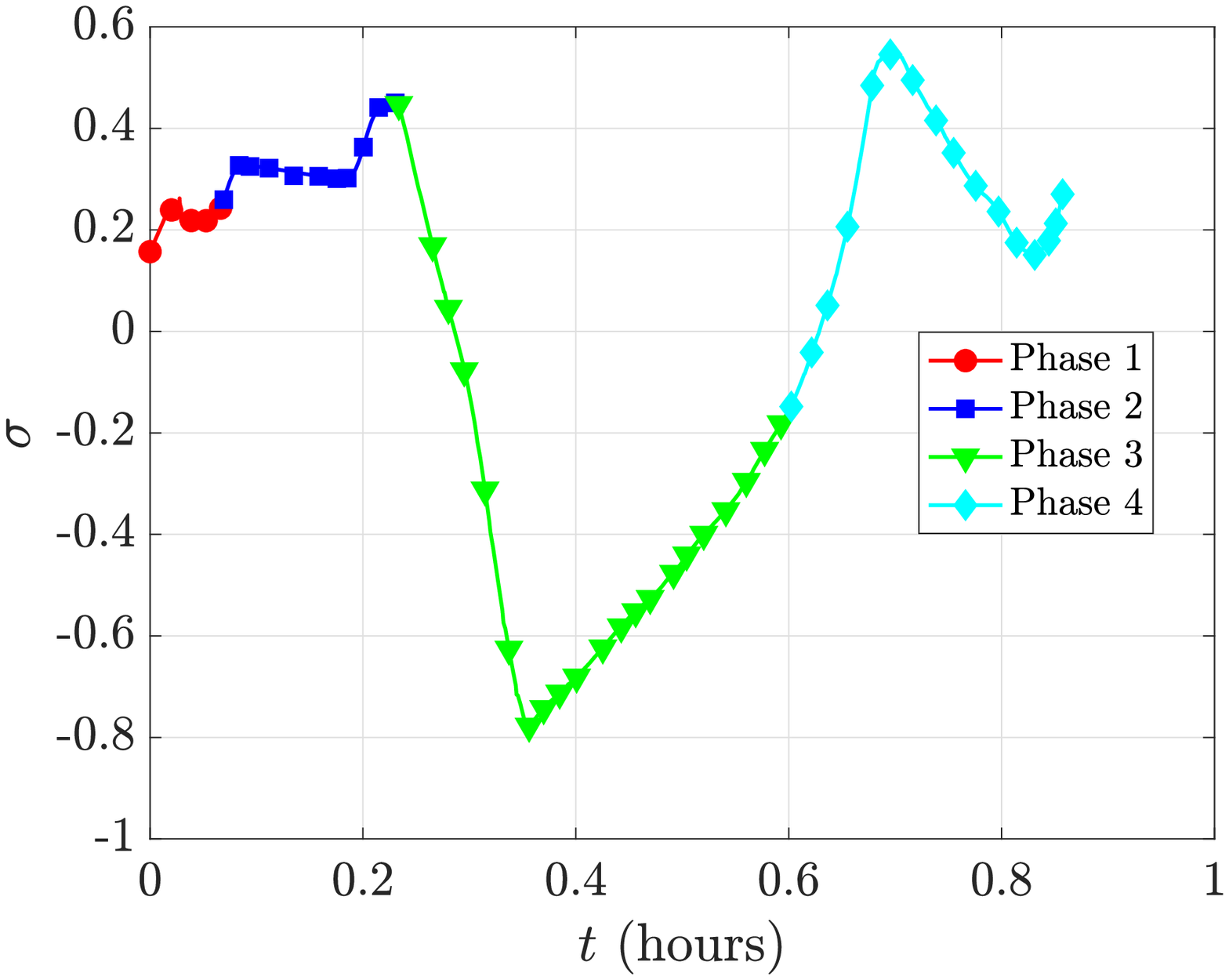}} \\
\subfloat[Chance Constrained]{\includegraphics[height = 2.1in]{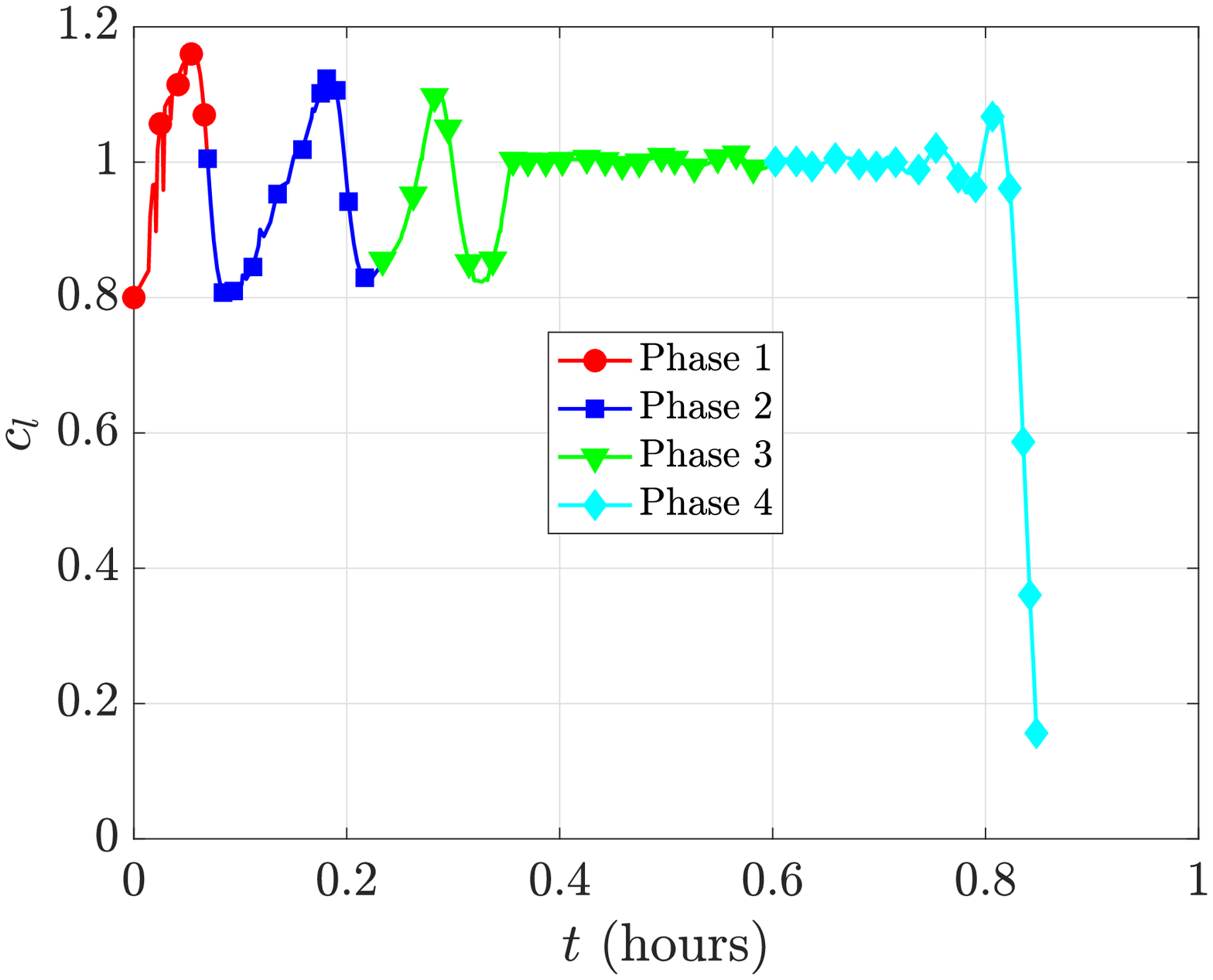}} 
~~~~\subfloat[Deterministic]{\includegraphics[height = 2.1in]{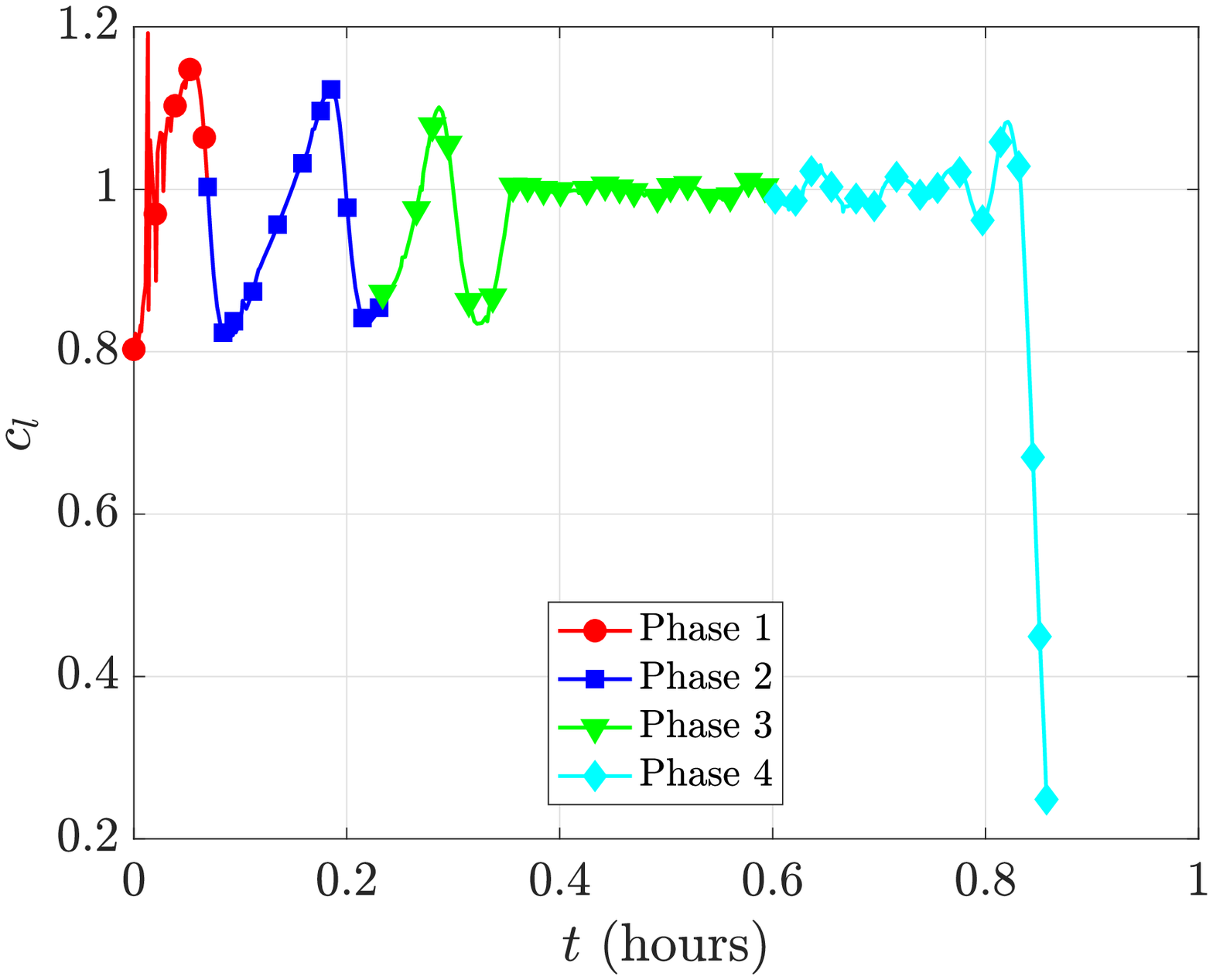}} \\
\caption{Controls for Example 2 and deterministic variation of Example 2.}
\label{fig:Jorris3DControls}
\end{figure}

This section provides results for solving Example 2 using the approach of Section~\ref{sect:review}, along with the warm start method from Section~\ref{sect:tech}.  The values  $\phi = 5\times 10^{-5}$ and $w = 1$ were used, and the Gaussian kernel was the starting kernel.  Additionally, it is noted that both the chance constrained and deterministic formulations of Example 2 were solved with the $\mathbb{GPOPS-II}$ setup discussed in Section~\ref{sect:setup} with no limit on the number of mesh refinement iterations required to reach a solution with a user-specified accuracy tolerance.  The solutions shown in Figs.~\ref{fig:Jorris3DStates} and~\ref {fig:Jorris3DControls} correspond to a single run using the Epanechnikov kernel and a single run when solving the deterministic version of this example.  The figures indicate that the solutions for the chance constrained and deterministic formulations are similar, and so the results obtained for the chance constrained formulation are reasonable.  It is further noted that the solutions shown in Figs.~\ref{fig:Jorris3DStates} and~\ref {fig:Jorris3DControls} indicate that, even though the formulation of Example 2 is similar to the formulation of Example 1, the solutions are quite different.  In particular, as shown in Fig.~\ref{fig:Jorris2DEpanech} for Example 1 and Fig.~\ref{fig:Jorris2DEpanech} for Example 2, the control for Example 1 has a bang-bang structure, while the controls for Example 2 are smooth.  This difference in the behavior of the control is why tuning of the starting bandwidth was required to obtain solutions for Example 1 ($w = 100$), but not for Example 2. 

Next, Table~\ref{table_ex2_final} compares the following results: (1) results obtained using the approach of Section~\ref{sect:review} for twenty runs of the chance constrained version of the example [Eqs.~\eqref{eq:costvers3}--\eqref{eq:ex3path} and~\eqref{eq:CC2_ex2}] applying two different kernels, and (2) results obtained for twenty runs of the deterministic formulation of the example [Eqs.~\eqref{eq:costvers3}--\eqref{eq:ex3path} and \eqref{eq:ex3path_nofly}].  The results in Table~\ref{table_ex2_final} indicate that the average optimal costs were lower for the chance constrained formulation than for the deterministic formulation of Example 2, for reasons discussed in Section~\ref{sect:highmesh}.  Moreover, as discussed in Section~\ref{sect:highmesh}, because of the use of sampling in the chance constrained formulation of Example 2, the run times are higher than the run times for the deterministic formulation.  

\begin{table}[ht]
\caption{Results for Example 2. \label{table_ex2_final}}
\renewcommand{\baselinestretch}{1}\small\normalfont
\centering
\begin{tabular}{| c || c | c | c | c | }  \hline
  &  Split-Bernstein & Epanechnikov & Deterministic \\ \hline
  $\mu_{J^*}$ & $3052.414$ s & $3052.407$ s & $ 3086.738$ s \\ \hline
  $ \sigma_{J^*}$  & $0.0252$ s & $0.0248$ s & $0$ s \\ \hline
  $\mu_T$  & $55.733$ s & $63.066$ s  & $9.424$ s \\ \hline
  $T_{\min}$ & $39.771$ s & $45.075$ s  & $9.210$ s \\ \hline 
  $T_{\max}$ & $80.866$ s & $81.973$ s  & $9.584$ s \\ \hline
\end{tabular}
\end{table} 

Now comparing the results for Example 1 from Table~\ref{table_ex1_final} to the results for Example 2 from Table~\ref{table_ex2_final}, the run times are slightly higher for Example 2.  This difference is due to Example 2 being a more complex problem than Example 1.  The results indicate that two complex CCOCPs were efficiently solved using the approach of Ref~\ref{sect:review} along with the warm start method developed in Section~\ref{sect:tech}.

\section{Discussion}\label{sect:discussion}

The results of Sections \ref{sect:discusstech} and \ref{sect:examples} demonstrate the capabilities of the warm start method developed in Section \ref{sect:tech}.  In particular, the warm start method developed in Section~\ref{sect:tech} was applied effectively to solve two complex CCOCPs given in Sections~\ref{sect:discusstech} and~\ref{sect:examples} using biased KDEs and LGR collocation.  Moreover, it was found in Section \ref{sect:discusstech} that solving Example 1 using the warm start method was far more reliable and computationally efficient than solving Example 1 without a warm start (Section~\ref{sect:Jorris2Dnaive}). 

Now, while the warm start method developed in this paper is found to improve reliability and computational efficiency when solving CCOCPs, it is important to note several aspects of the method that must be implemented carefully.  First, for the two components of the method that were described in Sections \ref{sect:bandwidth-tuning} and \ref{sect:kernel-switching}, it is important to choose an appropriate starting bandwidth and kernel.  In particular, choosing an inappropriate starting bandwidth and kernel can result in the NLP solver not converging.  Moreover, with an inappropriate choice of a starting bandwidth, the NLP solver may converge to a solution different from the optimal solution.  Also, tuning the bandwidth and kernel can be time consuming.  It is noted, however, that if the trial runs for determining an appropriate starting bandwidth and kernel are performed using a small sample set and with limited mesh refinement iterations, results can be obtained rather quickly.  Additionally, when tuning the bandwidth and kernel using trial runs, the process can be terminated as soon as the NLP solver is found to not converge to a solution for one of the runs.  Thus, the maximum number of trail runs is not used until after an appropriate bandwidth and kernel combination has been found.  It is also noted that convergence to an infeasible solution is always possible because a slightly different solution to the CCOCP is obtained with each run due to the use of a different sample set for each run.  As a result, obtaining an infeasible solution in the trial runs is probable.  The choice of starting kernel can, however, affect how often the NLP solver converges to an infeasible solution.  It is further noted that the NLP solver will sometimes shift from this infeasible solution to a feasible solution when the bandwidth and kernel are switched.  

Next, examining the third component of the method as described in Section \ref{sect:sample-size-increasing}, the size of the starting sample set can affect whether or not a solution is obtained.  When the initial sample size is too small, important features of the samples such as modes, mean, and range can be lost.  Thus, when the starting sample size is too small, the key features for the smaller sample set will be different from those of the larger sample sets.  As a result, the solution of the NLP obtained using larger sample sizes may have different properties from the solution of the NLP obtained using the smaller sample size.  This difference can, in turn, lead to the NLP solver not converging to a solution when the number of samples is increased.  Additionally, it was found that, even if the NLP solver converges, more computation effort may be required on the mesh refinement iteration when the sample size is first increased.  Finally, there was large variation in the amount of time required for the mesh refinement iteration where the sample set is switched to the full sample.  Increasing the sample size more gradually (that is, by using more increments with smaller changes between increments) can potentially decrease the time for the mesh refinement iteration where the full sample size is used.  Conversely, by using a greater number of increments, the total computation time may start to increase, even if the time for that one mesh refinement iteration is reduced.  Additionally, this computation time may be increased even further if the greater number of increments results in extra mesh refinement iterations (as discussed in Section~\ref{sect:sample-size-increasing}).

\section{Conclusions}\label{sect:conclude}

A warm start method has been developed to increase the efficiency of solving chance constrained optimal control problems using biased kernel density estimators and Legendre-Gauss-Radau collocation.  First, through a motivating example, it was shown that solving a chance constrained optimal control problem without a warm start can be unreliable and computationally inefficient.  Using the computational issues of solving this example as a starting point, the warm start method has been developed.  The warm start method consists of three components that are designed to aid convergence of the NLP solver, while simultaneously decreasing the required computation time and reducing sensitivity to the kernel and the initial guess.  These three components of the warm start method are: bandwidth tuning, kernel switching, and incremental sample size increasing.  The warm start method has then been applied to solve the motivating chance constrained optimal control problem using biased kernel density estimators and Legendre-Gauss-Radau collocation.  Finally, a second and more complex variation of this chance constrained optimal control problem has also been solved with the warm start method, and the results analyzed.  The results show that the warm start method developed in this paper has the potential to significantly improve reliability and computational efficiency when solving complex chance constrained optimal control problems.

\section*{Acknowledgments}

The authors gratefully acknowledge support for this research from the from the U.S.~National Science Foundation under grants CMMI-1563225, DMS-1522629, and DMS-1819002.  

\renewcommand{\baselinestretch}{1.0}
\normalsize\normalfont
\bibliographystyle{aiaa}


\begin{thebibliography}{10}
\newcommand{\enquote}[1]{``#1''}

\bibitem{Betts3}
Betts, J.~T., {\em Practical Methods for Optimal Control and Estimation Using
  Nonlinear Programming\/}, SIAM Press, Philadelphia, 2nd ed., 2009.

\bibitem{Keil1}
Keil, R., Aggarwal, R., Kumar, M., and Rao, A.~V., \enquote{Application of
  Chance-Constrained Optimal Control to Optimal Obstacle Avoidance,} {\em AIAA
  Guidance, Navigation and Control Conference, San Diego\/}, January 2019,
  pp.~\url{https://doi.org/10.2514/6.2019--0647}.

\bibitem{Kumar1}
Zhao, Z. and Kumar, M., \enquote{Split-Bernstein Approach to Chance-Constrained
  Optimal Control,} {\em Journal of Guidance, Control, and Dynamics\/},
  Vol.~40, No.~11, November 2017, pp.~2782--2795.
  {\url{https://doi.org/10.2514/1.G002551}}.

\bibitem{Caillau}
Caillau, J.-B., Cerf, M., Sassi, A., Trelat, E., and Zidani, H.,
  \enquote{Solving chance constrained optimal control problems in aerospace via
  Kernel Density Estimation,} {\em Optimal Control Applications and Methods\/},
  Vol.~39, No.~5, Wiley 2018, pp.~1833--1858.
  {\url{https://doi.org/10.1002/oca.2445}}.

\bibitem{blackmore10}
Blackmore, L., Ono, M., Bektassov, A., and Williams, B.~C., \enquote{A
  Probabilistic Particle-Control Approximation of Chance-Constrained Stochastic
  Predictive Control,} {\em IEEE Transactions on Robotics\/}, Vol.~26, No.~3,
  June 2010, pp.~502--517. {\url{https://doi.org/10.1109/TRO.2010.2044948}}.

\bibitem{Blackmore1}
Blackmore, L., Ono, M., and Williams, B.~C., \enquote{Chance-Constrained
  Optimal Path Planning With Obstacles,} {\em IEEE Transactions on Robotics\/},
  Vol.~27, No.~6, December 2011, pp.~1080--1094.
  {\url{https://doi.org/10.1109/TRO.2011.2161160}}.

\bibitem{ono10}
Ono, M., Blackmore, L., and Williams, B.~C., \enquote{Chance Constrained Finite
  Horizon Optimal Control with Nonconvex Constraints,} {\em Proceedings of the
  2010 {{American Control Conference}}\/}, {IEEE}, {Baltimore, MD}, 2010, pp.
  1145--1152. {\url{https://doi.org/10.1109/ACC.2010.5530976}}.

\bibitem{okamoto19}
Okamoto, K. and Tsiotras, P., \enquote{Optimal Stochastic Vehicle Path Planning
  Using Covariance Steering,} {\em IEEE Robotics and Automation Letters\/},
  Vol.~4, No.~3, July 2019, pp.~2276--2281.
  {\url{https://arxiv.org/abs/1809.03380}}.

\bibitem{hokayem13}
Hokayem, P., D.~Chatterjee, D., and Lygeros, J., \enquote{Chance-constrained
  LQG with bounded control policies,} {\em 52nd IEEE Conference on Decision and
  Control\/}, Florence, Italy, December 2013, pp. 2471--2476.
  {\url{https://doi.org/10.1109/CDC.2013.6760251}}.

\bibitem{Pinter}
Pint\'{e}r, J., \enquote{Deterministic Approximations of Probability
  Inequalities,} {\em Zeitschrift f\"{u}r Operations Research\/}, Vol.~33,
  No.~4, July 1989, pp.~219--239. {\url{https://doi.org/10.1007/BF01423332}}.

\bibitem{muhlpfordt18}
Muhlpfordt, T., Faulwasser, T., and Hagenmeyer, V., \enquote{A generalized
  framework for chance-constrained optimal power flow,} {\em Sustainable
  Energy, Grids and Networks\/}, Vol.~16, 2018, pp.~231--242.
  {\url{https://doi.org/10.1016/j.segan.2018.08.002}}.

\bibitem{Nemirovski}
Nemirovski, A. and Shapiro, A., \enquote{Convex Approximations of Chance
  Constrained Programs,} {\em SIAM Journal on Optimization\/}, Vol.~17, No.~4,
  November 2006, pp.~969--996. {\url{https://doi.org/10.1137/050622328}}.

\bibitem{Pagnoncelli1}
Pagnoncelli, B.~K., Ahmed, S., and Shapiro, A., \enquote{Sample Average
  Approximation Method for Chance Constrained Programming: Theory and
  Applications,} {\em Journal of Optimization Theory and Application\/},
  Vol.~142, 2009, pp.~399--416. {\url{
  https://doi.org/10.1007/s10957--009--9523--6}}.

\bibitem{ono15}
Ono, M., Pavone, M., Kuwata, Y., and Balaram, J., \enquote{Chance-Constrained
  Dynamic Programming with Application to Risk-Aware Robotic Space
  Exploration,} {\em Autonomous Robots\/}, Vol.~39, No.~4, Dec. 2015,
  pp.~555--571. {\url{https://doi.org/10.1007/s10514--015--9467--7}}.

\bibitem{Calafiore1}
Calafiore, G.~C. and Campi, M.~C., \enquote{The Scenario Approach to Robust
  Control Design,} {\em IEEE Transactions on Automatic Control\/}, Vol.~51,
  No.~5, May 2006, pp.~742--753.
  {\url{https://doi.org/10.1109/TAC.2006.875041}}.

\bibitem{Calafiore2}
Calafiore, G.~C. and Fagiano, L., \enquote{Robust Model Predictive Control via
  Scenario Optimization,} {\em IEEE Transactions on Automatic Control\/},
  Vol.~58, No.~1, January 2013, pp.~219--224.
  {\url{https://doi.org/10.1109/TAC.2012.22}}.

\bibitem{Campi1}
Campi, M.~C. and Garatti, S.~A., \enquote{A Sampling-and-Discarding Approach to
  Chance-Constrained Optimization: Feasibility and Optimality,} {\em Journal of
  Optimization Theory and Applications\/}, Vol.~148, 2011, pp.~257--280.
  {\url{https://doi.org/10.1007/s10957--010--9754--6}}.

\bibitem{Chai}
Chai, R., Savvaris, A., Tsuordos, A., Chai, S., Xia, Y., and Wang, S.,
  \enquote{Solving Trajectory Optimization Problems in the Presence of
  Probabilistic Constraints,} {\em IEEE Transactions on Cybernetics\/},
  February 2019, pp.~1--14. {\url{https://doi.org/10.1109/TCYB.2019.2895305}}.

\bibitem{Ahmed}
Ahmed, S., \enquote{Convex relaxations of chance constrained optimization
  problems,} {\em Optimization Letters\/}, Vol.~8, No.~1, January 2014,
  pp.~1--12. {\url{https://doi.org/10.1007/s11590--013--0624--7}}.

\bibitem{Calfa}
Calfa, B.~A., Grossman, I.~E., Agarwal, A., Bury, S.~J., and Wassick, J.~M.,
  \enquote{Data-driven individual and joint chance-constrained optimization via
  kernel smoothing,} {\em Computers and Chemical Engineering\/}, Vol.~78, July
  2015, pp.~51--69. {\url{https://doi.org/10.1016/j.compchemeng.2015.04.012}}.

\bibitem{Keil2}
Keil, R., Miller, A.~T., Kumar, M., and Rao, A.~V., \enquote{Method for Solving
  Chance Constrained Optimal Control Problems Using Biased Kernel Density
  Estimators,} {\em arXiv\/}, March 2020,
  pp.~\url{https://arxiv.org/abs/2003.08010}.

\bibitem{Benson2}
Benson, D.~A., Huntington, G.~T., Thorvaldsen, T.~P., and Rao, A.~V.,
  \enquote{{D}irect {T}rajectory {O}ptimization and {C}ostate {E}stimation via
  an {O}rthogonal {C}ollocation {M}ethod,} {\em Journal of Guidance, Control,
  and Dynamics\/}, Vol.~29, No.~6, November-December 2006, pp.~1435--1440.
  {\url{https://doi.org/10.2514/1.20478}}.

\bibitem{Rao8}
Rao, A.~V., Benson, D.~A., Darby, C.~L., Francolin, C., Patterson, M.~A.,
  Sanders, I., and Huntington, G.~T., \enquote{{A}lgorithm 902: {GPOPS}, {A}
  {MATLAB} {S}oftware for {S}olving {M}ultiple-{P}hase {O}ptimal {C}ontrol
  {P}roblems {U}sing the {G}auss {P}seudospectral {M}ethod,} {\em ACM
  Transactions on Mathematical Software\/}, Vol.~37, No.~2, April--June 2010,
  Article 22, 39 pages. {\url{https://doi.org/10.1145/1731022.1731032}}.

\bibitem{Garg1}
Garg, D., Patterson, M.~A., Darby, C.~L., Francolin, C., Huntington, G.~T.,
  Hager, W.~W., and Rao, A.~V., \enquote{{D}irect {T}rajectory {O}ptimization
  and {C}ostate {E}stimation of {F}inite-{H}orizon and {I}nfinite-{H}orizon
  {O}ptimal {C}ontrol {P}roblems via a {R}adau {P}seudospectral {M}ethod,} {\em
  Computational Optimization and Applications\/}, Vol.~49, No.~2, June 2011,
  pp.~335--358. {\url{https://doi.org/10.1007/s10589--009--9291--0}}.

\bibitem{Garg2}
Garg, D., Patterson, M.~A., Hager, W.~W., Rao, A.~V., Benson, D.~A., and
  Huntington, G.~T., \enquote{A Unified Framework for the Numerical Solution of
  Optimal Control Problems Using Pseudospectral Methods,} {\em Automatica\/},
  Vol.~46, No.~11, November 2010, pp.~1843--1851.
  {\url{https://doi.org/10.1016/j.automatica.2010.06.048}}.

\bibitem{Patterson2015}
Patterson, M.~A., Hager, W.~W., and Rao, A.~V., \enquote{A $ph$ Mesh Refinement
  Method for Optimal Control,} {\em Optimal Control Applications and
  Methods\/}, Vol.~36, No.~4, July--August 2015, pp.~398--421.
  {\url{https://doi.org/10.1002/oca.2114}}.

\bibitem{HagerHouRao15a}
Hager, W.~W., Hou, H., and Rao, A.~V., \enquote{Lebesgue Constants Arising in a
  Class of Collocation Methods,} {\em IMA Journal of Numerical Analysis\/},
  Vol.~13, No.~1, October 2017, pp.~1884--1901.
  \url{https://doi.org/10.1093/imanum/drw060}.

\bibitem{HagerHouRao16a}
Hager, W.~W., Hou, H., and Rao, A.~V., \enquote{Convergence Rate for a Gauss
  Collocation Method Applied to Unconstrained Optimal Control,} {\em Journal of
  Optimization Theory and Applications\/}, Vol.~169, No.~3, 2016, pp.~801--824.
  {\url{https://doi.org/10.1007/s10957--016--0929--7}}.

\bibitem{HagerLiuMohapatraWangRao19}
Hager, W.~W., Liu, J., Mohapatra, S., Rao, A.~V., and Wang, X.-S.,
  \enquote{Convergence Rate for a Gauss Collocation Method Applied to
  Constrained Optimal Control,} {\em SIAM Journal on Control and
  Optimization\/}, Vol.~56, No.~2, 2018, pp.~1386--1411.
  \url{https://doi.org/10.1137/16M1096761}.

\bibitem{DuChenHager2019}
Du, W., Chen, W., Yang, L., and Hager, W.~W., \enquote{Bounds for Integration
  Matrices That Arise in Gauss and Radau Collocation,} {\em Computational
  Optimization and Applications\/}, Vol.~74, September 2019, pp.~259--273.
  {\url{https://doi.org/10.1007/s10589--019--00099--5}}.

\bibitem{Kumar3}
{Zhao}, Z. and {Kumar}, M., \enquote{A split-bernstein approach to chance
  constrained programs,} {\em 53rd IEEE Conference on Decision and Control\/},
  2014, pp. 6621--6626. {\url{https://doi.org/10.1109/CDC.2014.7040428}}.

\bibitem{Kumar2}
{Zhao}, Z. and {Kumar}, M., \enquote{A MCMC/Bernstein approach to chance
  constrained programs,} {\em 2014 American Control Conference\/}, 2014, pp.
  4318--4323. {\url{https://doi.org/10.1109/ACC.2014.6859159}}.

\bibitem{MCMCMethods}
Neal, R.~M., {\em Handbook of Markov Chain Monte Carlo\/}, chap.~5, CRC Press,
  Boca Raton, Florida, 2011, pp. 113--162.

\bibitem{Gill1}
Gill, P.~E., Murray, W., and Saunders, M.~A., \enquote{{SNOPT}: {A}n {SQP}
  {A}lgorithm for {L}arge-{S}cale {C}onstrained {O}ptimization,} {\em SIAM
  Review\/}, Vol.~47, No.~1, 2005, pp.~99--131.
  {\url{https://doi.org/10.1137/S0036144504446096}}.

\bibitem{Gill2}
Gill, P.~E., Wong, E., Murray, W., and Saunders, M.~A., {\em User's Guide for
  SNOPT Version 7.6: Software for Large Scale Nonlinear Programming\/}, January
  2017. {\url{http://www.sbsi-sol-optimize.com/manuals/SNOPT%20Manual.pdf}}.

\bibitem{Biegler2}
Biegler, L.~T. and Zavala, V.~M., \enquote{{L}arge-{S}cale {N}onlinear
  {P}rogramming {U}sing {IPOPT}: {A}n {I}ntegrating {F}ramework for
  {E}nterprise-{W}ide {O}ptimization,} {\em Computers and Chemical
  Engineering\/}, Vol.~33, No.~3, March 2008, pp.~575--582.
  \url{https://doi.org/10.1016/j.compchemeng.2008.08.006}.

\bibitem{Byrd1}
Byrd, R.~H., Nocedal, J., and Waltz, R.~A., \enquote{KNITRO: An Integrated
  Package for Nonlinear Optimization,} {\em Large Scale Nonlinear
  Optimization\/}, Springer Verlag, 2006, pp. 35--59.
  {\url{https://doi.org/10.1007/0--387--30065--1_4}}.

\bibitem{Jorris2}
Jorris, T.~R. and Cobb, R.~G., \enquote{Multiple Method 2-D Trajectory
  Optimization Satisfying Waypoints and No-Fly Zone Constraints,} {\em Journal
  of Guidance, Control, and Dynamics\/}, Vol.~31, No.~3, May--June 2008,
  pp.~543--553. {\url{https://doi.org/10.2514/1.32354}}.

\bibitem{Garg3}
Garg, D., Hager, W.~W., and Rao, A.~V., \enquote{Pseudospectral Methods for
  Solving Infinite-Horizon Optimal Control Problems,} {\em Automatica\/},
  Vol.~47, No.~4, April 2011, pp.~829--837.
  {\url{https://doi.org/10.1016/j.automatica.2011.01.085}}.

\bibitem{Liu2015}
Liu, F., Hager, W.~W., and Rao, A.~V., \enquote{{A}daptive {M}esh {R}efinement
  for {O}ptimal {C}ontrol {U}sing {N}onsmoothness {D}etection and {M}esh {S}ize
  {R}eduction,} {\em Journal of the Franklin Institute\/}, Vol.~352, No.~10,
  October 2015, pp.~4081--4106.
  {\url{https://doi.org/10.1016/j.jfranklin.2015.05.028}}.

\bibitem{Liu2018}
Liu, F., Hager, W.~W., and Rao, A.~V., \enquote{{A}daptive {M}esh {R}efinement
  for {O}ptimal {C}ontrol {U}sing {D}ecay {R}ates of {L}egendre {P}olynomial
  {C}oefficients,} {\em IEEE Transactions on Control System Technology\/},
  Vol.~26, No.~4, 2018, pp.~1475--1483.
  {\url{https://doi.org/10.1109/TCST.2017.2702122}}.

\bibitem{Darby2}
Darby, C.~L., Hager, W.~W., and Rao, A.~V., \enquote{{A}n $hp$-{A}daptive
  {P}seudospectral {M}ethod for {S}olving {O}ptimal {C}ontrol {P}roblems,} {\em
  Optimal Control Applications and Methods\/}, Vol.~32, No.~4, July--August
  2011, pp.~476--502. {\url{ https://doi.org/10.1002/oca.957}}.

\bibitem{Darby3}
Darby, C.~L., Hager, W.~W., and Rao, A.~V., \enquote{{D}irect {T}rajectory
  {O}ptimization {U}sing a {V}ariable {L}ow-{O}rder {A}daptive {P}seudospectral
  {M}ethod,} {\em Journal of Spacecraft and Rockets\/}, Vol.~48, No.~3,
  May--June 2011, pp.~433--445. {\url{https://doi.org/10.2514/1.52136}}.

\bibitem{Francolin2014a}
Francolin, C.~C., Hager, W.~W., and Rao, A.~V., \enquote{{C}ostate
  {A}pproximation in {O}ptimal {C}ontrol {U}sing {I}ntegral {G}aussian
  {Q}uadrature {C}ollocation {M}ethods,} {\em Optimal Control Applications and
  Methods\/}, Vol.~36, No.~4, July--August 2015, pp.~381--397. {\url{
  https://doi.org/10.1002/oca.2112}}.

\bibitem{Patterson2014}
Patterson, M.~A. and Rao, A.~V., \enquote{$\mathbb{GPOPS-II}$, {A} {MATLAB}
  {S}oftware for {S}olving {M}ultiple-{P}hase {O}ptimal {C}ontrol {P}roblems
  {U}sing $hp$-{A}daptive {G}aussian {Q}uadrature {C}ollocation {M}ethods and
  {S}parse {N}onlinear {P}rogramming,} {\em ACM Transactions on Mathematical
  Software\/}, Vol.~41, No.~1, October 2014, pp.~1--37.
  {\url{https://doi.org/10.1145/2558904}}.

\bibitem{Patterson2012}
Patterson, M.~A. and Rao, A.~V., \enquote{{E}xploiting {S}parsity in {D}irect
  {C}ollocation {P}seudospectral {M}ethods for {S}olving {C}ontinuous-{T}ime
  {O}ptimal {C}ontrol {P}roblems,} {\em Journal of Spacecraft and Rockets,\/},
  Vol.~49, No.~2, March--April 2012, pp.~354--377.
  {\url{https://doi.org/10.2514/1.A32071}}.

\bibitem{bowman1}
Bowman, A. and Azzalini, A., {\em Applied Smoothing Techniques for Data
  Analysis: The Kernel Approach with S-Plus Illustrations\/}, Oxford
  Statistical Science Series, OUP Oxford, 1997.

\bibitem{silverman1}
Silverman, B.~W., {\em Density Estimation for Statistics and Data Analysis\/},
  Monographs on Statistics and Applied Probability 26, Chapman and Hall/CRC
  Press, 1986.

\bibitem{Neal2}
Neal, R.~M., {\em Probabilistic Inference Using Markov Chain Monte Carlo
  Methods\/}, Department of Computer Science, University of Toronto, September
  1993, Technical Report, CRG-TR-93-1.
  {\url{https://bayes.wustl.edu/Manual/RadfordNeal.review.pdf}}.

\bibitem{Epanech1}
Epanechnikov, V.~A., \enquote{Nonparametric Estimation of a Multivariate
  Probability Density,} {\em Theory of Probability and Its Applications\/},
  Vol.~14, No.~1, 1969, pp.~156--161. {\url{https://doi.org/10.1137/1114019}}.

\bibitem{Roycet1}
Roycet, J.~O. and Polak, E., \enquote{Implementable Algorithm for Stochastic
  Optimization Using Sample Average Approximations,} {\em Journal of
  Optimization Theory and Application\/}, Vol.~122, No.~1, July 2004,
  pp.~157--184. {\url{https://doi.org/10.1023/B:JOTA.0000041734.06199.71}}.

\bibitem{Jorris3}
Jorris, T.~R. and Cobb, R.~G., \enquote{Three-Dimensional Trajectory
  Optimization Satisfying Waypoint and No-Fly Zone Constraints,} {\em Journal
  of Guidance, Control, and Dynamics\/}, Vol.~32, No.~2, March--April 2009,
  pp.~551--572. {\url{https://doi.org/10.2514/1.37030}}.

\end{thebibliography}

\end{document}